%% file: main.tex
\pdfoutput=1
\documentclass[toc]{mathprint}
\usepackage{rutarmacros}
\usepackage{macros}
\title[Multifractal analysis via Lagrange duality]{Multifractal analysis via Lagrange duality}

\author{Alex Rutar}{
    Department of Mathematics,
    The University of British Columbia,
    1984 Mathematics Road,
    Vancouver BC V6T 1Z2,
    Canada
}{alex@rutar.org}

\begin{document}
\begin{abstract}
    We provide a self-contained exposition of the well-known multifractal formalism for self-similar measures satisfying the strong separation condition.
    At the heart of our method lies a pair of quasiconvex optimization problems which encode the parametric geometry of the Lagrange dual associated with the constrained variational principle.
    We also give a direct derivation of the Hausdorff dimension of the level sets of the upper and lower local dimensions by exploiting certain weak uniformity properties of the space of Bernoulli measures.
\end{abstract}

\section{Introduction}

\subsection{Prelude on multifractal analysis}
Many naturally-occurring systems can be modelled by a sequence of independent and identically distributed random variables, or more generally by a dynamical system supporting an invariant measure.
In these situations, one is often interested in the long-term behaviour given \emph{typical} initial conditions.

In \emph{multifractal analysis}, one is interested in the \emph{non-typical} behaviour of a dynamical system.
One of the oldest examples which hints at multifractal analysis can be found in the study of \emph{Besicovitch--Eggleston sets} \cite{zbl:0009.05301,zbl:0031.20801,zbl:0014.25702}.
Fix an integer $b\geq 2$ and consider the base-$b$ expansion of a point $x\in[0,1]$.
Suppose $x$ has base-$b$ expansion $x=0.a_1a_2a_3\ldots$ where $a_i\in\{0,1,\ldots,b-1\}$.
Then the \emph{digit frequencies} of $x$ exist almost surely: there is a set $E\subset[0,1]$ of full Lebesgue measure such that for each $x=0.a_1a_2a_3\ldots\in E$ and $i=0,1,\ldots,b-1$,
\begin{equation*}
    \lim_{n\to\infty}\frac{\#\{j:a_j = i,1\leq i\leq n\}}{n}=\frac{1}{b}.
\end{equation*}
But what can be said about non-typical behaviour---the complement of $E$ is non-empty, but what can be said about its structure?

In order to understand non-typical behaviour, it is reasonable to associate with the original system different invariant measures, for which typical points of the new measure will reveal new properties of the underlying system.
For instance, we may define a probability measure on $[0,1]$ by taking points with base-$b$ expansion chosen randomly according to some alternative distribution $(p_0,\ldots,p_{b-1})$ where, independently for each $n$, we choose $a_n$ to be the digit $i$ with probability $p_i$.
Typical points for this measure will not be elements of the set $E$ (unless $p_0=\cdots=p_{b-1}=1/b$), but the measure still has positive dimension: in fact, the set of typical points for this measure has Hausdorff dimension precisely
\begin{equation*}
    \frac{-\sum_{i=0}^{b-1} p_i\log p_i}{\log b}.
\end{equation*}
In particular, even though the complement of $E$ has measure $0$, it has Hausdorff dimension 1.
Moreover, this gives us a good understanding of the size of the sets of points where the digit frequencies exist and are given by some distribution.
But beyond this, there are still many points (in fact also a set of full Hausdorff dimension) for which the digit frequencies do not exist at all and one can ask a variety of questions about such points.

In this paper, we will focus on the multifractal analysis of measures, of which the above consideration is a special case.
Given a compactly supported Borel probability measure in $\R^d$ and $x\in\supp\mu$, we wish to study the local scaling properties of the measure $\mu$ at $x$ from the perspective of the \emph{local dimension}
\begin{equation*}
    \dim_{\loc}(\mu,x)=\lim_{r\to 0}\frac{\log\mu(B(x,r))}{\log r},
\end{equation*}
when the limit exists (or taking a limit infimum or limit supremum when it does not).
If the measure $\mu$ is sufficiently nice, then $\dim_{\loc}(\mu,x)$ almost surely exists and attains some constant value.
Moreover, this average scaling property is intimately related to the \emph{dimension} of the measure $\mu$, since a probability measure which has larger local dimensions is more concentrated in neighbourhoods of points in its support.
But what about non-typical points?
A reasonable question to ask is how large such sets are: we will be most interested in the quantity
\begin{equation}\label{e:multifractal-def}
    f_\mu(\alpha)=\dim\{x\in\supp\mu:\dim_{\loc}(\mu,x)=\alpha\}
\end{equation}
where $\dim$ is either Hausdorff or packing dimension; or perhaps some variant by considering only the points where the lower or upper limit converge to some value $\alpha$.
We refer to the function $f_\mu$ as the \emph{multifractal spectrum} of the measure $\mu$.
For more detail on the Hausdorff and packing dimensions, we refer to reader to the books \cite{zbl:1390.28012,zbl:0869.28003}.
The level sets of local dimensions will often be dense subsets of the support of the measure, so other notions of dimension such as the box dimensions (which are stable under closure) do not typically yield meaningful information.

Note that the function $f_\mu$ appears to be defined in a somewhat pathological way: it is the Hausdorff dimension of a level set of an asymptotically-defined pointwise function, so there appears to be no \emph{a priori} reason to expect that the dependence on $\alpha$ is sensible.
Nonetheless, the \emph{multifractal miracle} (perhaps first observed non-rigorously in the physics literature \cite{zbl:1184.37028}) is that the quantity $f_\mu$, for the class of measures that we will consider, is in fact a non-trivial concave analytic function of $\alpha$.
We will study the multifractal quantity $f_\mu$, in conjunction with other quantitative measures of smoothness such as the $L^q$-spectrum which will be introduced in later sections.

For digit frequency measures, and more generally for the family of measures that we will consider, we emphasize that the multifractal formalism has been rigorously established in \cite{zbl:0763.58018,zbl:0873.28003}.
Expositions of the standard proofs can be found in \cite[§11.2]{zbl:0869.28003} and \cite[Chapter~5]{zbl:1543.28001}.
Our goal in this paper is to provide a novel approach to and exposition of these results in order to motivate the multifractal analysis of measures and highlight the relationship with Lagrange duality, constrained optimization, and large deviations theory.

\subsection{Constrained optimization and duality}
Recall our general goal: we have fixed a measure $\mu$ and we wish to understand the multifractal spectrum $f_\mu(\alpha)$ defined in \cref{e:multifractal-def}.
Suppose that we could choose a measure $\nu$ such that $\dim_{\loc}(\mu,x)=\alpha$ for $\nu$-a.e.\ $x\in\supp\nu$.
Then we would immediately obtain the lower bound
\begin{equation*}
    f_\mu(\alpha)\geq\dimH\mu\coloneqq\sup\{\dimH E:E\subset\supp\mu:\mu(E)=1\}
\end{equation*}
(see \cref{ss:dimensions} for details on dimensions of measures).

Unfortunately, this idea is not particularly helpful in general, since on any set with given Hausdorff dimension one can find measures with dimension arbitrarily close to the dimension of the original set.
However, measures defined on digit sets (and more generally the measures that we will consider in this paper) have more structure: they are projections of \emph{invariant} measures on some compact metric space $\Omega$ under a coding map $\pi\colon\Omega\to\R^d$.
Given a continuous transformation $T\colon \Omega\to \Omega$, recall that a measure $\nu$ on $\Omega$ is \emph{$T$-invariant} if for all Borel sets $E\subset \Omega$, $\nu(E)=\nu(T^{-1}(E))$, and \emph{ergodic} if $\nu(E)\in\{0,1\}$ for any set $E$ with $E=T^{-1}(E)$.
Let $\mathcal{M}$ denote the space of $T$-invariant and ergodic measures.
Then one might ask for the following \emph{constrained variational principle}: is the general inequality
\begin{equation}\label{e:var-principle}
    f_\mu(\alpha)\geq\sup_{\nu\in\mathcal{M}}\left\{\dimH \pi_*\nu:\dim_{\loc}(\mu,x)=\alpha\text{ for $\pi_*\nu$-a.e.\ }x\right\}
\end{equation}
sharp?
Here, $\pi_*\nu$ is the pushforward measure of $\nu$ under the projection map $\pi$.
Moreover, the \emph{dimension} of $\pi_*\nu$ is determined in a loose sense by the \emph{entropy} of $\nu$ as well as the geometry of the projection map $\pi$.
If such a variational principle holds, the situation seems more hopeful!
In fact, assuming that the maps $T$ and $\pi$ are reasonably well-behaved, our domain and objective function satisfies important forms of regularity: the space $\mathcal{M}$ is the set of extreme points of a compact and convex metric space (with topology precisely the weak-star topology) and the entropy map on $\mathcal{M}$ is upper semicontinuous \cite[Proposition~5.2 and Theorem~8.2]{zbl:0475.28009}.
We can therefore replace the general question of understanding the multifractal spectrum with two new questions:
\begin{enumerate}[nl]
    \item Can we understand the variational lower bound in \cref{e:var-principle}; and
    \item Can we prove that equality in \cref{e:var-principle} actually holds?
\end{enumerate}

Associated with the constrained optimization problem which gives a lower bound to $f_\mu(\alpha)$ is the corresponding \emph{unconstrained Lagrange dual} minimization problem defining $\tau_\mu(q)$, which encodes the Lagrange multipliers problem associated with the constrained maximization.
These dual optimization problems are fundamentally related by the concave conjugate.
When the underlying space and the objective function are particularly well-behaved, the relationship between $f_\mu(\alpha)$ and $\tau_\mu(q)$ is well-understood (a particularly elegant case occurs when the domain and objective function are convex, see for instance \cite[Part~III]{zbl:0193.18401}).
However, the key difference in our perspective is to make minimal assumptions on the domain and objective functions, and to instead understand the \emph{parametric geometry} of the maximization problem \emph{as a function of $\alpha$}.
This allows us to convert topological properties of the Lagrange dual or generic properties of the function $\tau_\mu$ to information about the primal function $f_\mu(\alpha)$ (for instance, connectedness of the set of minimizers or differentiability of $\tau_\mu$ at a value $q$ implies the concave conjugate relation at values of $\alpha$ corresponding to $q$).
This approach is introduced in an abstract setting in \cref{ss:optim}.

We then need to prove that the inequality in \cref{e:var-principle} is in fact an equality.
To do this, we make a key observation: in our setting, the unconstrained dual corresponds precisely to the \emph{$L^q$-spectrum} of the measure (see \cref{ss:lq-spectrum} for a formal definition).
This proof can be found in \cref{t:lq-form}, which uses the \emph{method of types} from large deviations theory to establish a variational formula for the function $\tau_\mu(q)$.
The author was motivated to take this approach by the recent work of Kolossváry \cite{zbl:1549.37013}.
In order to complete the proof, it then remains to establish a general concave conjugate relationship between the $L^q$-spectrum and the multifractal spectrum through a \emph{geometric large deviations bound} (see \cref{p:lq-concave-bound}).
This part of the argument has much more of a geometric flavour, since we must handle points which are ``invisible'' to the dynamics (in that they are non-typical for any pushforward of an ergodic measure on $\Omega$).
In \cref{p:var-upper} we also give a direct argument using a variant of a classical density theorem for Hausdorff dimension to establish equality in \cref{e:var-principle}.

In fact, the method of types argument and the density argument are precisely the tools required to reduce a covering argument for general points to a covering argument for ``typical'' points, which allows us to take advantage of the underlying dynamical structure of the sets and measures under consideration.
Such reductions are necessary since covering arguments must work for all points of a set simultaneously.

In this paper, we focus on the multifractal analysis of the local dimensions of a probability measure $\mu$.
One can ask analogous versions of this question, such as for Birkhoff averages on a dynamical system.
In this case, the constrained and unconstrained optimization problems considered above are intimately related to the notion of \emph{topological pressure} (see \cite{zbl:1319.37016} for a survey on this dynamical perspective to multifractal analysis).
The general approach to duality is certainly very applicable in this dynamical context as well.
In fractal geometry, however, a key role is played by the projection map $\pi$ and the geometry of Euclidean space which prevents us from taking a purely dynamical perspective.

Beyond the scope of this article, there are many new difficulties but also many recent breakthroughs, such as the establishment of a \emph{dynamical dimension gap} for self-affine sets \cite{zbl:1387.37026} (in some sense the strongest possible failure of equality in \cref{e:var-principle}) or the formula for the $L^q$-spectrum of an overlapping self-similar measure \cite{zbl:1426.11079} for $q>1$, which draws on deep ideas from discretized sum-product theory and additive combinatorics.

Many of the ideas in these notes can be extended and used in more general settings.
Some results (such as the general setup in \cref{ss:optim}, or the approach to the variational principle in \cref{ss:alt}) have been stated with this use-case in mind.
However, generally speaking, for clarity of exposition and in order to highlight the many connections between these diverse areas of mathematics while minimizing the technical complexity, we will restrict our attention to a special setting: the case of self-similar measures satisfying the strong separation condition.

\subsection{Self-similar sets and measures}
We now introduce the class of measures which we will focus on for the remainder of this paper.
This class of measures generalizes the ``digit frequency measures'' discussed in the introduction.
A key distinction, which only causes minor difficulties in our setting but is a critical difficulty in developing this theory in a more general context, is that we do not require the subdivision into parts to all be of the same size.

Fix a finite index set $\mathcal{I}$.
For each $i\in\mathcal{I}$, let $S_i\colon\R^d\to\R^d$ be a \defn{contracting similarity}, that is for each $i\in\mathcal{I}$ there is some $r_i\in(0,1)$ so that
\begin{equation*}
    |S_i(x)-S_i(y)|=r_i\cdot|x-y|\text{ for all }x,y\in\R^d.
\end{equation*}
We then call the finite family $\{S_i\}_{i\in\mathcal{I}}$ an \defn{iterated function system of similarities}, or \defn{IFS} for short.

Let $\mathcal{P}=\mathcal{P}(\mathcal{I})\subset \R^{\mathcal{I}}$ denote the set of probability vectors
\begin{equation*}
    \mathcal{P}=\bigl\{(p_i)_{i\in\mathcal{I}}:p_i\in[0,1],\textstyle{\sum_{i\in\mathcal{I}}p_i=1}\bigr\}.
\end{equation*}
We equip the space $\mathcal{P}$ with the metric from $\R^{\mathcal{I}}$: note that $\mathcal{P}$ is compact.
Now to each self-similar IFS, there is a unique compact set $K$ satisfying
\begin{equation*}
    K=\bigcup_{i\in\mathcal{I}}S_i(K)
\end{equation*}
and, given a probability vector $\bm{p}=(p_i)_{i\in\mathcal{I}}\in\mathcal{P}$, a unique Borel probability measure $\mu_{\bm{p}}$ satisfying
\begin{equation*}
    \mu_{\bm{p}}=\sum_{i\in\mathcal{I}}p_i\cdot\mu_{\bm{p}}\circ S_i^{-1}.
\end{equation*}
We call the set $K$ the \defn{attractor} and the measure $\mu_{\bm{p}}$ the \defn{invariant measure} associated with the IFS $\{S_i\}_{i\in\mathcal{I}}$ and probabilities $(p_i)_{i\in\mathcal{I}}$.
\begin{remark}
    For convenience, only for the reference measure $\mu_{\bm{p}}$, we will assume without loss of generality that $p_i>0$ for every $i\in\mathcal{I}$, so that $\supp\mu_{\bm{p}}=K$.
\end{remark}

A classic example of an iterated function system is the Cantor IFS, which is the system $\{x\mapsto x/3,x\mapsto x/3+2/3\}$.
Then the attractor $K$ is the middle-thirds Cantor set (i.e.\ the set of points in $[0,1]$ with ternary expansion consisting of only $0$s and $2$s).
If we take the probability vector $(1/2,1/2)\in\mathcal{P}$, then the corresponding measure $\mu$ is just normalized Hausdorff $\frac{\log 2}{\log 3}$-measure restricted to $K$ and normalized to have measure $1$.

We say that the IFS $\{S_i\}_{i\in\mathcal{I}}$ satisfies the \defn{strong separation condition}, or \defn{SSC} for short, if $S_i(K)\cap S_j(K)=\varnothing$ for all $i\neq j$ in $\mathcal{I}$.
Note that, formally, the SSC disallows measures which arise from a digit frequency set as discussed in the introduction.
However, those measures satisfy the following slightly weaker condition.
We say that the IFS $\{S_i\}_{i\in\mathcal{I}}$ satisfies the \defn{open set condition}, or \defn{OSC} for short, if there is an open set $U$ such that $S_i(U)\subset U$ for all $i\in\mathcal{I}$, and $S_i(U)\cap S_j(U)=\varnothing$ for $i\neq j$.

The main results in this paper also hold under the OSC, though there are some additional technical complications related to non-injectivity of the projection map $\pi$ that we can bypass under the SSC.
For simplicity of exposition, we will for the remainder of the document assume that the SSC holds.

\subsection{Shift space and symbolic coding}
As a result of the iterative nature of the construction of self-similar sets and measures, it is natural to introduce notation formalizing this iterative procedure.
One natural setting for such a construction is the \emph{full shift} on symbols $\mathcal{I}$.
Let $\Omega=\mathcal{I}^{\N}$ denote the set of sequences on symbols $\mathcal{I}$ equipped with the product topology, and let $\pi\colon\Omega\to K$ denote the natural coding map defined by
\begin{equation*}
    \bigl\{\pi\bigl((i_n)_{n=1}^\infty\bigr)\bigr\} = \lim_{n\to\infty} S_{i_1}\circ\cdots\circ S_{i_n}(K).
\end{equation*}
Note that if the IFS satisfies the SSC, the map $\pi\colon\Omega\to K$ is a homeomorphism.

It will be useful to introduce some notation to handle sequences in $\Omega$.
First, let $\mathcal{I}^*=\bigcup_{n=0}^\infty\mathcal{I}^n$ denote the space of all finite words on the alphabet $\mathcal{I}$, equipped with the operation of concatenation.
Given $\mtt{i}\in\mathcal{I}^*$ and $\mtt{k}\in\mathcal{I}^*$, we say that $\mtt{i}$ is a \emph{prefix} of $\mtt{k}$ if there is a $\mtt{j}\in\mathcal{I}^*$ so that $\mtt{k}=\mtt{i}\mtt{j}$.
Concatenation also extends naturally to expressions of the form $\mtt{i}\gamma$ for $\gamma\in\Omega$, and allows us to define the notion of a prefix of an infinite word.
In particular, for $\gamma\in\Omega$, we let $\gamma\npre{n}$ denote the unique prefix of $\gamma$ of length $n$.
As usual, for $\mtt{i} = (i_1,\ldots,i_n)\in\mathcal{I}^*$, we write
\begin{align*}
    r_{\mtt{i}} &= r_{i_1}\cdots r_{i_n},\\
    p_{\mtt{i}} &= p_{i_1}\cdots p_{i_n},\\
    S_{\mtt{i}} &= S_{i_1}\circ\cdots\circ S_{i_n}.
\end{align*}

We equip the sequence space $\Omega$ with the \emph{shift map} $\sigma\colon\Omega\to\Omega$ defined by
\begin{equation*}
    \sigma(i_1,i_2,i_3,\ldots)=(i_2,i_3,\ldots).
\end{equation*}
For example, for all $n\in\N$, $\gamma=\gamma\npre{n}\sigma^n(\gamma)$.
A probability vector $\bm{p}\in\mathcal{P}$ then naturally defines a shift-invariant and ergodic product measure $\bm{p}^{\N}$.
We call such a measure \defn{Bernoulli}.
The self-similar measure $\mu_{\bm{p}}$ is simply the pushforward of $\bm{p}^{\N}$ by the coding map $\pi$.

Finally, we set for $\bm{p}=(p_i)_{i\in\mathcal{I}}$
\begin{equation*}
    \Omega_{\bm{p}}\coloneqq\Bigl\{(i_n)_{n=1}^\infty\in\Omega:\lim_{n\to\infty}\frac{\#\{\ell:i_\ell=j\text{ for }1\leq \ell\leq n\}}{n}=p_j\text{ for }j\in\mathcal{I}\Bigr\},
\end{equation*}
in other words the collection of sequences in $\Omega$ where the digit frequencies exist and are given by the probability vector $\bm{p}$.

It will turn out, for our purposes, that it suffices to consider only Bernoulli measures, and so we will not consider the more general class of ergodic measures on $\Omega$ for the remainder of the document.
Note that for more general iterated function systems (with corresponding geometric complexities associated with the map $\pi$), this is no longer the case.

\section{Dimensions and \texorpdfstring{$L^q$}{Lq}-spectra of self-similar measures}
In this section, we begin with some classical theory: the notion of the Hausdorff dimension of a measure.
We will also derive a variational formula for the $L^q$-spectrum of a self-similar measure satisfying the SSC using the \emph{method of types} from large deviations theory.
Such a strategy was motivated by the main results and techniques used in \cite{zbl:1549.37013}.

There are substantial advantages to deriving a formula for the $L^q$-spectrum in terms of a variational formula, in contrast to classical proofs such as the one given in \cite[§11.2]{zbl:0869.28003} which work directly with an analytic closed formula.
Firstly, it suffices to consider continuous functions rather than smooth functions, which gives much more flexibility in more complex settings.
Secondly, as we will see in \cref{ss:optim} and \cref{ss:param}, there are many indirect techniques to understand the geometry of a general optimization problem, which makes it easier to obtain precise information instead of directly using an implicit formula.
Such techniques become essential in the settings where a closed formula for the $L^q$-spectrum may not even exist.

\subsection{Hausdorff dimension of self-similar measures}\label{ss:dimensions}
Given a Borel probability measure $\mu$, we denote the \defn{local dimension} of $\mu$ at $x$ by
\begin{equation*}
    \dim_{\loc}(\mu,x)=\lim_{r\to 0}\frac{\log\mu(B(x,r))}{\log r}
\end{equation*}
when the limit exists.
If the limit does not exist, we write $\underline{\dim}_{\loc}(\mu,x)$ (resp.\ $\overline{\dim}_{\loc}(\mu,x)$) to denote the respective quantity except with a limit infimum (resp.\ limit supremum) in place of the limit.
We can use the local dimensions to give a formulation of the Hausdorff dimension of $\mu$.
\begin{definition}
    Given a compactly supported Borel probability measure $\mu$, we define
    \begin{equation*}
        \dimH\mu=\esssup_{x\in\supp\mu}\underline{\dim}_{\loc}(\mu,x).
    \end{equation*}
    We say that $\mu$ is \defn{exact-dimensional} if there is an $\alpha\in\R$ so that $\dim_{\loc}(\mu,x)=\alpha$ for $\mu$-a.e.\ $x\in\supp\mu$.
\end{definition}
Equivalently, as explained in \cite[Proposition~10.1]{zbl:0869.28003},
\begin{equation*}
    \dimH\mu=\sup\{\dimH E:E\subset\supp\mu:\mu(E)=1\}.
\end{equation*}

Now given a self-similar IFS $\{S_i\}_{i\in\mathcal{I}}$ and $\bm{p}=(p_i)_{i\in\mathcal{I}},\bm{w}=(w_i)_{i\in\mathcal{I}}\in\mathcal{P}(\mathcal{I})$, denote the \defn{cross entropy} and \defn{entropy} by
\begin{equation*}
    H(\bm{w},\bm{p})=\sum_{i\in\mathcal{I}}w_i\log(1/p_i)\qquad\text{and}\qquad H(\bm{w})=H(\bm{w},\bm{w}).
\end{equation*}
We also denote the \defn{Kullback--Leibler divergence} (also known as the \emph{relative entropy}) by
\begin{equation*}
    \DKL{\bm{w}}{\bm{p}}=H(\bm{w},\bm{p})-H(\bm{w}).
\end{equation*}
We recall by Jensen's inequality applied to the logarithm, $\DKL{\bm{w}}{\bm{p}}\geq 0$ with equality if and only if $\bm{w}=\bm{p}$.
For a general introduction to these concepts in information theory, we refer the reader to \cite{zbl:1140.94001}.
Finally, we define the \defn{Lyapunov exponent}
\begin{equation*}
    \chi(\bm{w})=\sum_{i\in\mathcal{I}}w_i\log(1/r_i).
\end{equation*}
The Lyapunov exponent captures the asymptotic contraction rate at points typical for the measure $\mu_{\bm{w}}$.
The key role of the Lyapunov exponent in dimension theory is highlighted, for instance, by the seminal works of Ledrappier--Young \cite{zbl:0605.58028,zbl:1371.37012}.

We will now determine the dimensions of self-similar measures.
We first observe the following technical lemma, which allows a reduction to balls which intersect precisely one image $S_\mtt{i}(K)$ for $\mtt{i}\in\mathcal{I}^*$.
In essence, this will allow us to treat the geometry of $\mu_{\bm{p}}$ in a purely symbolic way.
Note that $\diam(K)\cdot r_\mtt{i}=\diam(S_\mtt{i}(K))$.
\begin{lemma}\label{l:symb-reduction}
    For all sufficiently small $\delta>0$, there is a constant $c=c(\bm{p},\delta)>0$ so that for all $\mtt{i}\in\mathcal{I}^*$ and $x\in S_\mtt{i}(K)$,
    \begin{equation*}
        c\cdot p_{\mtt{i}}\leq\mu_{\bm{p}}\bigl(B(x,\delta\cdot r_{\mtt{i}})\bigr)\leq p_{\mtt{i}}.
    \end{equation*}
\end{lemma}
\begin{proof}
    Since the IFS satisfies the strong separation condition, for all sufficiently small $\delta>0$ the $\delta\cdot r_\mtt{i}$-neighbourhood of $S_\mtt{i}(K)$ in $K$ is again $S_\mtt{i}(K)$.
    Thus
    \begin{equation*}
        \mu_{\bm{p}}\bigl(B(x,\delta\cdot r_{\mtt{i}})\bigr)\leq p_{\mtt{i}}.
    \end{equation*}
    On the other hand, for each fixed $\delta>0$, there is an $N\in\N$ (depending on $\delta$) and a word $\mtt{j}\in\mathcal{I}^N$ so that
    \begin{equation*}
        S_{\mtt{i}}\circ S_{\mtt{j}}(K)\subseteq B(x,\diam(S_\mtt{i}(K))\cdot\delta).
    \end{equation*}
    Taking $c=\min\{p_{\mtt{j}}:\mtt{j}\in\mathcal{I}^N\}$ gives the desired result.
\end{proof}
We also obtain a simple lemma which gives information about local dimensions in terms of the digit frequencies of the symbolic representation.
\begin{lemma}\label{l:symb-formula}
    Let $\bm{w}\in\mathcal{P}(\mathcal{I})$ and let $\gamma=(i_n)_{n=1}^\infty\in\Omega_{\bm{w}}$.
    Then
    \begin{equation*}
        \dim_{\loc}(\mu_{\bm{p}},\pi(\gamma)) = \frac{H(\bm{w},\bm{p})}{\chi(\bm{w})}.
    \end{equation*}
\end{lemma}
\begin{proof}
    By \cref{l:symb-reduction}, there is a $\delta>0$ and a $c>0$ so that for any $n\in\N$,
    \begin{equation*}
        c\cdot p_{\gamma\npre{n}}\leq \mu_{\bm{p}}\bigl(B(\pi(\gamma),\delta\cdot r_{\gamma\npre{n}})\bigr)\leq p_{\gamma\npre{n}}.
    \end{equation*}
    Let $\bm{q}_n$ denote the digit frequencies of $\gamma\npre{n}$, so that
    \begin{equation*}
        \dim_{\loc}(\mu_{\bm{p}},\pi(\gamma)) = \lim_{n\to\infty}\frac{\log p_{\gamma\npre{n}}}{\log r_{\gamma\npre{n}}} = \lim_{n\to\infty}\frac{H(\bm{q}_n,\bm{p})}{\chi(\bm{q}_n)} = \frac{H(\bm{w},\bm{p})}{\chi(\bm{w})}
    \end{equation*}
    from the definition of $\Omega_{\bm{w}}$, as claimed.
\end{proof}
In particular, we obtain the following well-known dimension formula for self-similar measures.
\begin{proposition}\label{p:rel-dim}
    Let $\{S_i\}_{i\in\mathcal{I}}$ be an IFS satisfying the SSC and let $\bm{p},\bm{w}\in\mathcal{P}$.
    Then for $\mu_{\bm{w}}$-a.e.\ $x\in\supp\mu_{\bm{w}}$,
    \begin{equation*}
        \dim_{\loc}(\mu_{\bm{p}},x)=\frac{H(\bm{w},\bm{p})}{\chi(\bm{w})}.
    \end{equation*}
    In particular, $\mu_{\bm{p}}$ is exact-dimensional with dimension
    \begin{equation*}
        \dimH\mu_{\bm{p}}=\frac{H(\bm{p})}{\chi(\bm{p})}.
    \end{equation*}
\end{proposition}
\begin{proof}
    By Kolmogorov's Strong Law of Large Numbers, $\bm{w}^{\N}(\Omega_{\bm{w}})=1$.
    Thus the dimensional result follows from \cref{l:symb-formula}.
\end{proof}

\subsection{\texorpdfstring{$L^q$}{Lq}-spectra of self-similar measures}\label{ss:lq-spectrum}
Let $\mu$ be a compactly supported Borel probability measure and let $q\in\R$.
We write for $r\in(0,1)$
\begin{equation*}
    G_\mu(r,q)=\sup\left\{\,\sum_i\mu\bigl(B(x_i,r)\bigr)^q:\{x_i\}_i\text{ is a $2r$-separated subset of }\supp\mu\,\right\}.
\end{equation*}
We then denote the \defn{$L^q$-spectrum} of $\mu$ at $q$ by
\begin{equation*}
    \tau_\mu(q)=\liminf_{r\to 0}\frac{\log G_\mu(r,q)}{\log r}.
\end{equation*}

We recall the following standard result.
\begin{lemma}
    Let $\mu$ be a compactly supported Borel probability measure.
    Then $\tau_\mu(q)$ is a concave and increasing function of $q$.
    Moreover, $\tau_\mu(0)=-\dimB(\supp\mu)$ and $\tau_\mu(1)=0$.
\end{lemma}
\begin{proof}
    Let $q_1<q_2$.
    The observation that $\tau_\mu(q)$ is increasing follows since $\mu(B(x_i,r))^{q_1}\geq\mu(B(x_i,r))^{q_2}$ for any $x_i\in \supp\mu$.

    In addition, concavity is a standard application of Hölder's inequality: let $0<\lambda<1$, and then with Hölder's inequality with exponents $1/\lambda$ and $1/(1-\lambda)$,
    \begin{equation}\label{e:theta-ineq}
        \sum_i\mu\bigl(B(x_i,r)\bigr)^{\lambda q_1+(1-\lambda)q_2}\leq\left(\sum_i\mu\bigl(B(x_i,r)\bigr)^{q_1}\right)^\lambda\left(\sum_i\mu\bigl(B(x_i,r)\bigr)^{q_2}\right)^{1-\lambda}
    \end{equation}
    and taking suprema and logarithms yields concavity.

    Finally, using the equivalence of packings and coverings, we see that $\tau_\mu(0)=-\dimB(\supp\mu)$.
    Moreover, since $\mu$ is a probability measure there is a constant $c>0$ so that $0<c\leq G_\mu(r,1)\leq 1$ for all $r>0$.
    Indeed, to see that such a constant $c$ exists, first observe that there is a constant $m$ depending only on the ambient dimension $d$ and a cover $\mathcal{B}$ of $\supp\mu$ such that each ball intersects at most $m$ other balls.
    By a greedy algorithm (take the ball $B\in\mathcal{B}$ with maximal measure, discard at most $m$ other balls which intersect it, and repeat), we obtain a sub-collection $\mathcal{B}'\subset\mathcal{B}$ such that
    \begin{equation*}
        1 \leq \sum_{B\in\mathcal{B}}\mu(B) \leq (m+1)\sum_{B\in\mathcal{B}_j}\mu(B).
    \end{equation*}
    Thus we can take $c = 1/(m+1)$.
    In particular, it follows that $\tau_\mu(1)=0$.
\end{proof}

In our special setting, we can rewrite the $L^q$-spectrum essentially as a symbolic sum.
For notational simplicity, we will write $\tau_{\bm{p}}$ in place of $\tau_{\mu_{\bm{p}}}$, where the dependency on the underlying IFS is implicit.
Since $(r_{\gamma\npre{n}})_{n=1}^\infty$ is a strictly decreasing sequence for each $\gamma\in\Omega$, setting
\begin{equation*}
    \Lambda_r\coloneqq\{\mtt{i}\in\mathcal{I}^*\setminus\{\varnothing\}:r_\mtt{i}\leq r<r_{\mtt{i}^-}\},
\end{equation*}
we have that $\Lambda_r$ is a section of $\mathcal{I}^*$ for each $r>0$ (that is, for all $n$ sufficiently large and $\mtt{i}\in\mathcal{I}^n$, $\mtt{i}$ has a unique prefix in $\Lambda_r$).
Here, for non-empty $\mtt{i} = (i_1,\ldots,i_n)$, we write $\mtt{i}^- = (i_1,\ldots,i_{n-1})$.
Heuristically, $\Lambda_r$ is a symbolic representation of the collection of images $S_\mtt{i}(K)$ with diameter approximately $r$.

We now have the following immediate application of \cref{l:symb-reduction}.
\begin{lemma}\label{l:symb-sum-formula}
    Let $\bm{p}\in\mathcal{P}$ be arbitrary and $q\in\R$.
    Then there are constants $c_1,c_2>0$ (depending on $q$) and a constant $\delta>0$ so that for all $r\in(0,1)$
    \begin{equation*}
        c_1 G_{\mu_{\bm{p}}}(\delta r,q)\leq\sum_{\mtt{i}\in\Lambda_r}p_\mtt{i}^q\leq c_2 G_{\mu_{\bm{p}}}(\delta r,q).
    \end{equation*}
    In particular,
    \begin{equation*}
        \tau_{\bm{p}}(q)=\liminf_{r\to 0}\frac{\log\sum_{\mtt{i}\in\Lambda_r}p_\mtt{i}^q}{\log r}.
    \end{equation*}
\end{lemma}
\begin{proof}
    By \cref{l:symb-reduction}, there are constants $\delta>0$ and $c>0$ so that for each $r\in(0,1)$ and $\mtt{i}\in\Lambda_r$, there is some $x_\mtt{i}\in S_\mtt{i}(K)$ so that
    \begin{equation}\label{e:bd1}
        c\cdot p_\mtt{i}\leq \mu\bigl(B(x_\mtt{i},\delta r)\bigr)\leq p_{\mtt{i}}
    \end{equation}
    and moreover $B(x,\delta r)\cap K\subseteq S_\mtt{i}(K)$, so that $\{x_\mtt{i}\}_{\mtt{i}\in\Lambda_r}$ is a $2\delta r$-separated subset of $K$.

    On the other hand, there is some $N\in\N$ so that if $\{x_i\}_i$ is any $2\delta r$-separated subset of $K$, then there is some unique $\mtt{i}_i\in\Lambda_r$ so that $B(x_i,\delta r)\cap K\subseteq S_{\mtt{i}_i}(K)$, and some $\mtt{j}_i\in\mathcal{I}^n$ so that $S_{\mtt{i}_i}\circ S_{\mtt{j}_i}(K)\subseteq B(x_i,\delta r)$.
    Thus with $c'=\min\{p_{\mtt{j}}:\mtt{j}\in\mathcal{I}^n\}$,
    \begin{equation}\label{e:bd2}
        c'p_{\mtt{i}_i}\leq\mu\bigl(B(x_i,\delta r)\bigr)\leq p_{\mtt{i}_i}.
    \end{equation}

    The desired inequalities now follow from \cref{e:bd1} and \cref{e:bd2}.
\end{proof}
In our next result, we will use the method of types from large deviations theory: we refer the reader to the book by Dembo \& Zeitouni \cite{zbl:1177.60035} for more background and detail.
To this end, for $n\in\N$, $i\in\mathcal{I}$, and $\mtt{i}=(i_1,\ldots,i_n)\in\mathcal{I}^n$, let $\xi_i(\mtt{i})$ denote the \emph{frequency of letter $i$ in $\mtt{i}$}, that is
\begin{equation}\label{e:xi-def}
    \xi_i(\mtt{i})=\frac{\#\{j:1\leq h\leq n,i_j=i\}}{n}.
\end{equation}
We refer to the probability vector $\bm{\xi}(\mtt{i})\coloneqq (\xi_i(\mtt{i}))_{i\in\mathcal{I}}\in\mathcal{P}(\mathcal{I})$ as the \emph{type} associated with $\mtt{i}$.

Now given an $r\in(0,1)$, we let
\begin{equation*}
    \mathcal{T}_r=\bigl\{\bm{\xi}(\mtt{i}):\mtt{i}\in\Lambda_r\bigr\}
\end{equation*}
denote the set of all types at scale $r$.
Conversely, given a $\bm{w}\in\mathcal{T}_r$, we let
\begin{equation*}
    \mathcal{C}_r(\bm{w})=\bigl\{\mtt{i}\in\Lambda_r:\bm{\xi}(\mtt{i})=\bm{w}\bigr\}.
\end{equation*}
The main point of introducing types is that the values $r_{\mtt{i}}$, $|\mtt{i}|$, and $p_{\mtt{i}}$ (for $\bm{p}\in\mathcal{P}$) are constant on each type class.
Moreover, there are few types, and we can determine the size of each type class up to a sub-exponential error.
Here, for notational simplicity, we use Landau's asymptotic notation, i.e.\ given real-valued functions $f,g$ defined on some domain $A$, we say $f=O(g)$ if there is a constant $C>0$ so that $|f(a)/g(a)|\leq C$ for all $a\in A$.
\begin{lemma}\label{l:types}
    Let $r\in(0,1)$ be arbitrary.
    Then
    \begin{equation}\label{e:few-types}
        \frac{\#\mathcal{T}_r}{\log r}=O\left(\frac{\log\log(1/r)}{\log(1/r)}\right).
    \end{equation}
    and, for fixed $\bm{w}\in\mathcal{T}_r$ with $\bm{w}=\bm{\xi}(\mtt{i})$ for $\mtt{i}\in\Lambda_r$,
    \begin{equation}\label{e:type-count}
        \frac{\#\mathcal{C}_r(\bm{w})}{|\mtt{i}|}=H(\bm{w})+O\left(\frac{\log\log(1/r)}{\log(1/r)}\right).
    \end{equation}
\end{lemma}
\begin{proof}
    First, since $r_i\in(0,1)$ for all $i\in\mathcal{I}$, for all $r\in(0,1)$ and $\mtt{i}\in\Lambda_r$, $|\mtt{i}|\asymp\log(1/r)$.

    Since the type $\bm{\xi}(\mtt{i})$ does not depend on the order of the letters in $\mtt{i}$, a direct count (see \cite[Lemma~2.1.2]{zbl:1177.60035} for more detail) gives for each $m\in\N$
    \begin{equation*}
        \#\{\bm{\xi}(\mtt{i}):\mtt{i}\in\mathcal{I}^m\}\leq (m+1)^{\#\mathcal{I}}.
    \end{equation*}
    Thus there is a fixed constant $M>0$ so that
    \begin{align*}
        \#\mathcal{T}_r &\leq \sum_{m=0}^{\lceil M\log(1/r)\rceil}\#\{\bm{\xi}(\mtt{i}):\mtt{i}\in\mathcal{I}^m\}\\
                        &\leq \bigl(M\log(1/r)+1\bigr)\cdot(M\log(1/r)+2)^{\#\mathcal{I}}.
    \end{align*}
    Taking logarithms and dividing through by $\log(1/r)$ gives \cref{e:few-types}.

    Finally, \cref{e:type-count} is precisely \cite[Lemma~2.1.8]{zbl:1177.60035}, again using the fact that $|\mtt{i}|\asymp\log(1/r)$.
    Alternatively, one can directly apply Stirling's formula---that $\log(n!)=n\log n - n +O(\log n)$---and use the observation that the number of distinct permutations of a word $\mtt{i}=(i_1,\ldots,i_n)$ with $k_i=n\cdot \bm{\xi}(\mtt{i})$ for $i\in\mathcal{I}$ is $n!\cdot\bigl(\prod_{i\in\mathcal{I}}k_i!\bigr)^{-1}$.
\end{proof}
We can now establish our formula for the $L^q$-spectrum of a self-similar measure.
\begin{theorem}\label{t:lq-form}
    Let $\bm{p}\in\mathcal{P}$.
    Then the limit defining $\tau_{\bm{p}}$ exists, and moreover
    \begin{equation*}
        \tau_{\bm{p}}(q)=\inf_{\bm{w}\in\mathcal{P}}\left\{\frac{q H(\bm{w},\bm{p})-H(\bm{w})}{\chi(\bm{w})}\right\}.
    \end{equation*}
\end{theorem}
\begin{proof}
    We fix $\bm{p}\in\mathcal{P}$ and $q\in\R$.
    Suppose $r\in(0,1)$, fix $\bm{w}\in\mathcal{T}_r$, and let $m$ denote the common value of $|\mtt{i}|$ for $\mtt{i}\in\mathcal{C}_r(\bm{w})$.
    Then for any $\mtt{i}\in \mathcal{C}_{r}(\bm{w})$,
    \begin{equation*}
        p_\mtt{i}^q=\prod_{i\in\mathcal{I}}p_i^{q\cdot n\cdot w_i}\qquad\text{and}\qquad r_\mtt{i}=\prod_{i\in\mathcal{I}}r_i^{n\cdot w_i}.
    \end{equation*}
    It therefore follows that
    \begin{equation*}
        \sum_{\mtt{i}\in\Lambda_r}p_\mtt{i}^q=\sum_{\bm{w}\in\mathcal{T}_{r}}\# \mathcal{C}_{r}(\bm{w})\prod_{i\in\mathcal{I}}p_i^{q\cdot n\cdot w_i},
    \end{equation*}
    and moreover, for fixed $\bm{w}\in\mathcal{T}_{r}$, by \cref{e:type-count} in \cref{l:types},
    \begin{equation}\label{e:opt}
        \frac{\log\# \mathcal{C}_{r}(\bm{w})\prod_{i\in\mathcal{I}}p_i^{q\cdot n\cdot w_i}}{\log r}=\frac{q H(\bm{w},\bm{p})-H(\bm{w})}{\chi(\bm{w})}+O\left(\frac{\log\log(1/r)}{\log(1/r)}\right).
    \end{equation}
    Thus by \cref{e:few-types}, with $\bm{w}^{(r)}\in\mathcal{T}_{r}$ chosen to minimize the quantity in \cref{e:opt},
    \begin{align*}
        \frac{\log\sum_{\mtt{i}\in\Lambda_r}p_\mtt{i}^q}{\log r}= \frac{ q H(\bm{w}^{(r)},\bm{p})-H(\bm{w}^{(r)})}{\chi(\bm{w}^{(r)})}+O\left(\frac{\log\log(1/r)}{\log(1/r)}\right).
    \end{align*}
    But the map $\bm{w}\mapsto\frac{q H(\bm{w},\bm{p})-H(\bm{w})}{\chi(\bm{w})}$ is continuous, so by compactness of $\mathcal{P}$ and since $\mathcal{T}_{r}$ becomes arbitrarily dense in $\mathcal{P}$ as $r$ goes to $0$, it follows that the limit defining $\tau_{\bm{p}}(q)$ in fact exists and is equal to the desired expression.
\end{proof}

\subsection{An explicit formula for the \texorpdfstring{$L^q$}{Lq}-spectrum}
We now have a formula for $\tau_{\bm{p}}(q)$ as the solution to an optimization problem, but it would be nice to have an explicit formula.
Let us work out a reasonable guess for the formula for $\tau_{\bm{p}}(q)$ using Lagrange multipliers, for now sweeping any technical issues under the rug.

The space of probability vectors $\mathcal{P}$ is a subset of the affine hyperplane $\{\bm{w}\in\R^{\mathcal{I}}:\norm{\bm{w}}=1\}$.
Write
\begin{equation}\label{e:phi-def}
    \phi(q,\bm{w})=q H(\bm{w},\bm{p})-H(\bm{w}).
\end{equation}
Suppose $\bm{w}$ is a minimizing vector for $\tau_{\bm{p}}(q)$.
Informally applying the method of Lagrange multipliers, such a vector $\bm{w}\in\mathcal{P}$ must satisfy
\begin{equation*}
    \frac{\partial}{\partial w_i}\frac{\phi(q,\bm{w})}{\chi(\bm{w})}=\lambda\frac{\partial}{\partial w_i}(\norm{\bm{w}}-1) = \lambda
\end{equation*}
for some $\lambda\in\R$.
In particular, for each $i\in\mathcal{I}$,
\begin{equation*}
    \chi(\bm{w})\cdot \frac{\partial}{\partial w_i} \phi(q,\bm{w}) - \phi(q,\bm{w})\cdot \frac{\partial}{\partial w_i} \chi(\bm{w})=\lambda\cdot\chi(\bm{w})^2
\end{equation*}
which, after rearranging and computing the derivatives, yields
\begin{equation}\label{e:lagrange-constraint}
    q\log(1/p_i)+\log w_i+1=\lambda\cdot\chi(\bm{w})+\tau_{\bm{p}}(q)\log(1/r_i).
\end{equation}
Multiplying by $w_i$ and summing over all $i\in\mathcal{I}$,
\begin{equation*}
    \phi(q,\bm{w})+1=\lambda\cdot\chi(\bm{w})+\chi(\bm{w})\tau_{\bm{p}}(q).
\end{equation*}
But $\bm{w}$ is a minimizing vector, so $\phi(q,\bm{w})/\chi(\bm{w})=\tau_{\bm{p}}(q)$ so that $\lambda=1/\chi(\bm{w})$.
Substituting this back into \cref{e:lagrange-constraint} and rearranging for $w_i$ yields
\begin{equation*}
    w_i=p_i^q r_i^{-\tau_{\bm{p}}(q)}.
\end{equation*}
Therefore, we would guess that $\tau_{\bm{p}}(q)$ must be the solution to the equation
\begin{equation*}
    \sum_{i\in\mathcal{I}}p_i^q r_i^{-\tau_{\bm{p}}(q)}=1.
\end{equation*}
and the minimization is attained at the probability vector $\bigl(p_i^q r_i^{-\tau_{\bm{p}}(q)}\bigr)_{i\in\mathcal{I}}$.

While the derivation above is certainly not rigorous (though it can be made rigorous without too much trouble), it helps to motivate the following formula for the $L^q$-spectrum which we will prove is the correct formula through a very simple information-theoretic argument.
First, define the function
\begin{equation*}
    \psi(q,t)=\sum_{i\in\mathcal{I}}p_i^q r_i^{-t}.
\end{equation*}
Note that for fixed $q$, $t\mapsto\psi(q,t)$ is a strictly decreasing function of $t$ with $\lim_{t\to\infty}\psi(q,t)=0$ and $\lim_{t\to\infty}\psi(q,t)=\infty$.
Thus there is a unique value $T(q)$ so that $\psi(q,T(q))=1$; and moreover, by the analytic implicit function theorem, $T(q)$ is an analytic function of $q$.
Finally, we define the vector
\begin{equation*}
    \bm{z}(q)=\bigl(p_i^q r_i^{-T(q)}\bigr)_{i\in\mathcal{I}}
\end{equation*}
which is indeed a probability vector by the definition of $T(q)$.
This was the original formula established for the $L^q$-spectrum by Arbeiter \& Patzschke \cite{zbl:0873.28003}.
\begin{proposition}[\cite{zbl:0873.28003}]\label{p:tau-diff}
    Fix $\bm{p}\in\mathcal{P}$ and $q\in\R$: then $\tau_{\bm{p}}(q)$ is the unique solution to
    \begin{equation}\label{e:implicit-lq}
        \sum_{i\in\mathcal{I}}p_i^q r_i^{-\tau_{\bm{p}}(q)}=1.
    \end{equation}
    In particular, $\tau_{\bm{p}}$ is analytic on $\R$.
\end{proposition}
\begin{proof}
    Observe for $\bm{w}\in\mathcal{P}$ that
    \begin{align*}
        \frac{ q H(\bm{w},\bm{p})-H(\bm{w})}{\chi(\bm{w})} = T(q)+\frac{\DKL{\bm{w}}{\bm{z}(q)}}{\chi(\bm{w})}.
    \end{align*}
    But $\DKL{\bm{w}}{\bm{z}(q)}\geq 0$ with equality if and only if $\bm{w}=\bm{z}(q)$, so the minimization is attained uniquely at the probability vector $\bm{z}(q)$ with value $T(q)$.
    Thus the result follows by \cref{t:lq-form}.
\end{proof}

\subsection{Asymptotes of the \texorpdfstring{$L^q$}{Lq}-spectra}\label{ss:lq-asymptotes}
\begin{figure}[t]
    \centering
    \input{figures/lq_plot_asymp}
    \caption{Plot of the $L^q$-spectrum of the measure $\mu_{\bm{p}}$, along with its asymptotes.}
    \label{f:lq}
\end{figure}
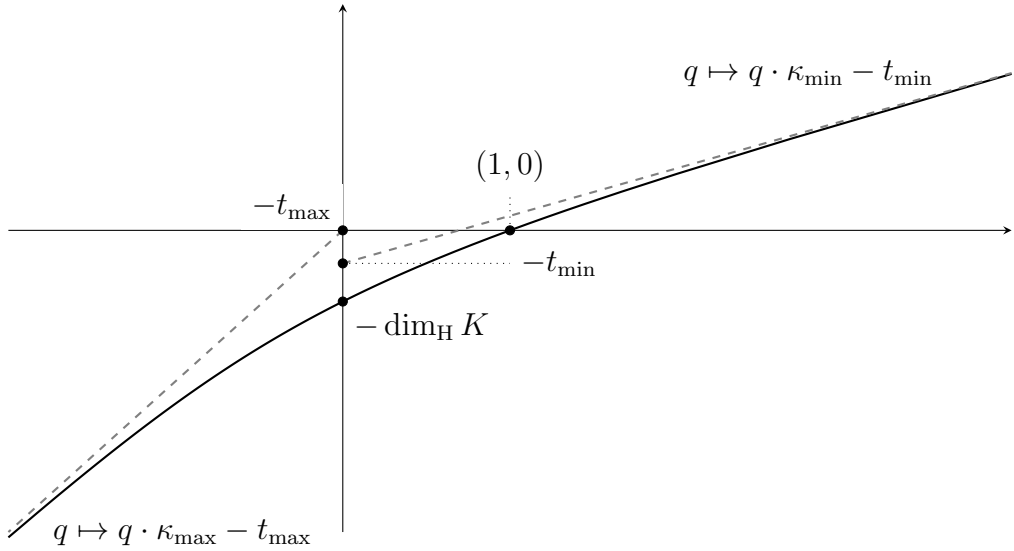

To conclude this section, let's analyze the asymptotic behaviour of $\tau_{\bm{p}}(q)$, as well as the limiting values of the corresponding optimization vector, which we recall is given by
\begin{equation*}
    \bm{z}(q)=\left(p_i^q r_i^{-\tau_{\bm{p}}(q)}\right)_{i\in\mathcal{I}}.
\end{equation*}
First, for each $i\in\mathcal{I}$, let $\kappa_i$ be such that $r_i^{\kappa_i}=p_i$, that is
\begin{equation*}
    \kappa_i = \frac{\log p_i}{\log r_i}.
\end{equation*}
In particular, $z_i(q)= r_i^{q\kappa_i -\tau_{\bm{p}}(q)}$.
But for $q>0$, $q\kappa_i-\tau_{\bm{p}}(q)\geq q\kappa_{\min}-\tau_{\bm{p}}(q)$, and the $z_i(q)$ cannot all converge to 0 as $q$ diverges to infinity.
The same argument holds for $q<0$, so that $\lim_{q\to\infty}z_i(q)>0$ if and only if $\kappa_i=\kappa_{\min}$, and $\lim_{q\to-\infty}z_i(q)>0$ if and only if $\kappa_i=\kappa_{\max}$.
Thus the limits $\bm{z}(\pm\infty)\coloneqq\lim_{q\to \pm\infty}\bm{z}(q)$ exist and are given by
\begin{equation*}
    z_i(-\infty)=\begin{cases}
        0 &: \kappa_i\neq\kappa_{\max}\\
        r_i^{s_{\max}} &: \kappa_i=\kappa_{\max}
    \end{cases}
    \qquad\text{and}\qquad
    z_i(\infty)=\begin{cases}
        0 &: \kappa_i\neq\kappa_{\min}\\
        r_i^{s_{\min}} &: \kappa_i=\kappa_{\min}
    \end{cases}
\end{equation*}
where
\begin{equation*}
    \sum_{\substack{i\in\mathcal{I}\\\kappa_i=\kappa_{\max}}}r_i^{s_{\max}}=1 \qquad\text{and}\qquad \sum_{\substack{i\in\mathcal{I}\\\kappa_i=\kappa_{\min}}}r_i^{s_{\min}}=1.
\end{equation*}
Equivalently,
\begin{equation*}
    s_{\max}=\lim_{q\to-\infty}(q\kappa_i - \tau_{\bm{p}}(q))
    \qquad\text{and}\qquad
    s_{\min}=\lim_{q\to\infty}(q\kappa_i - \tau_{\bm{p}}(q)).
\end{equation*}
This implies that $\tau_{\bm{p}}(q)$ has asymptotes given by
\begin{equation*}
    \ell_{-\infty}(q)=q\kappa_{\max} - s_{\max}\qquad\text{and}\qquad \ell_\infty(q)=q\kappa_{\min} -s_{\min}.
\end{equation*}
A general depiction of the function $\tau_{\bm{p}}(q)$ along with its asymptotes is given in \cref{f:lq}.

\section{Lagrange duality and multifractal formalism}\label{s:convex-duality}
Before proceeding with the remaining results of this section, let us summarize the current situation.
We are attempting to obtain a formula for the multifractal spectrum $f_{\bm{p}}(\alpha)$, for which we have the variational lower bound
\begin{equation*}
    f_{\bm{p}}(\alpha)\geq\sup_{\bm{w}\in\mathcal{P}}\left\{\dimH\mu_{\bm{w}}:\dim_{\loc}(\mu_{\bm{p}},x)=\alpha\text{ for $\mu_{\bm{p}}$-a.e.\ $x$}\right\}.
\end{equation*}
By \cref{p:rel-dim}, the optimization on the right hand side is precisely
\begin{equation}\label{e:constrained}
    \sup_{\bm{w}\in\mathcal{P}}\left\{\frac{H(\bm{w})}{\chi(\bm{w})}:\frac{H(\bm{w},\bm{p})}{\chi(\bm{w})}=\alpha\right\}.
\end{equation}
In order to understand this optimization problem, taking motivation from the method of Lagrange multipliers, we want to study the \emph{Lagrange dual}, which is optimization problem
\begin{equation}\label{e:dual}
    \inf_{\bm{w}\in\mathcal{P}}\left\{\frac{q\cdot H(\bm{w},\bm{p})-H(\bm{w})}{\chi(\bm{w})}\right\}.
\end{equation}
But this optimization problem is precisely the formula for the $L^q$-spectrum of $\mu_{\bm{p}}$, as proven in \cref{t:lq-form}!
Therefore, in order to complete the proof, we need two main ingredients:
\begin{enumerate}[nl]
    \item\label{im:duality} We need to understand the relationship between the dual optimization problems \cref{e:constrained} and \cref{e:dual}.
    \item\label{im:lq-bound} We need to obtain bounds on $f_{\bm{p}}(\alpha)$ by using the $L^q$-spectrum of $\mu_{\bm{p}}$.
\end{enumerate}
These two ingredients are handled in the subsequent sections.
In \cref{ss:optim}, we address the duality in \cref{im:duality} in a general context---real-valued continuous functions on a compact Hausdorff topological space.
Then in \cref{ss:Lq-bound} we study the relationship between $f_{\bm{p}}$ and $\tau_{\bm{p}}^*$: in fact, we will prove that $f_{\bm{p}}(\alpha)\leq\tau_{\bm{p}}^*(\alpha)$ in general.

Finally, we provide an alternative proof of the variational formula which avoids the machinery of duality in \cref{ss:alt}.
This proof is conceptually useful for different reasons, since it also immediately generalizes to different types of level sets of local dimensions.

\subsection{Continuous optimization and duality}\label{ss:optim}
Suppose $\Delta$ is a compact Hausdorff topological space and suppose we are given a continuous function $u\colon\Delta\to\R$ and an upper semicontinuous function $v\colon\Delta\to \R$.
We consider the constrained optimization
\begin{equation*}
    f(\alpha)=\max_{\bm{w}\in\Delta}\bigl\{v(\bm{w}):u(\bm{w})=\alpha\bigr\}
\end{equation*}
with corresponding unconstrained dual
\begin{equation}\label{e:tau-def}
    \tau(q)=\min_{\bm{w}\in\Delta}\bigl\{q\cdot u(\bm{w})-v(\bm{w})\bigr\}.
\end{equation}
Here, the maximum over the empty set is $-\infty$.
Note that both the maximum and minimum are attained on compact sets since $v$ is upper semicontinuous.
Of course, $\tau$ is concave since it is an infimum of affine functions.
On the other hand, $f$ need not be concave.

In some sense, one can think of the function $\tau(q)$ as encoding the geometry of the Lagrange multiplier problem associated with $f(\alpha)$.
However, this is only a motivating heuristic since in our abstract setup, there is no differentiable structure in sight.

Before we continue, let's recall some basic facts from convex optimization.
For a more in-depth introduction, we refer the reader to \cite{zbl:0193.18401}.
First, for a general function $g\colon\R\to\R\cup\{-\infty\}$, we denote the \emph{concave conjugate} by
\begin{equation*}
    g^*(\alpha) = \inf_{q\in\R}(q\alpha-g(q)).
\end{equation*}
Note that $g^*$ is always concave, and $g^{**}$ is the concave hull of $g$.

Now suppose moreover that $g\colon\R\to\R\cup\{-\infty\}$ is a concave function.
Then for $q\in\R$, we write $\partial g(q)$ to denote the \emph{subdifferential} of $g$ at $q$, i.e.
\begin{equation*}
    \partial g(q)=\{\alpha: \alpha(y-q)+g(q)\geq g(y)\text{ for any }y\in\R\}.
\end{equation*}
Equivalently,
\begin{equation}\label{e:subdiff-rel}
    g^*(\alpha)+g(q)\leq\alpha q
\end{equation}
with equality if and only if $\alpha\in\partial g(q)$.
We also let $\partial^-g(q)$ (resp.\ $\partial^+g(q)$) denote the left (resp.\ right) derivative of $g$ of $q$.
Then $\partial g(q)=[\partial^+g(q),\partial^-g(q)]$.
In particular, $g$ is differentiable at $q$ if and only if $\partial g(q)=\{\alpha\}$, in which case $g'(q)=\alpha$.
Since the subdifferentials form an ordered family of intervals of $\R$ which overlap only on their endpoints, there can be at most countably many points with non-singleton subdifferential: in particular, $g$ is differentiable at all but countably many $q\in\R$.

Now, we say that a line
\begin{equation*}
    \ell(t)=a\cdot t+b
\end{equation*}
is a \emph{supporting line} for $g$ at $q$ if $g(q)=\ell(q)$ and $\ell(t)\geq g(t)$ for all $t\in\R$.
Equivalently, $\partial g(q)$ is precisely the set of possible slopes of supporting lines for $g$ at $q$.

Now consider specifically the function $\tau$ from \cref{e:tau-def}.
We say that $\tau$ is \defn{supported} at $(q,\alpha)$ for $q\in\R$ if there is a $\bm{w}\in\Delta$ so that $\tau(q)=q\cdot u(\bm{w})-v(\bm{w})$ and $u(\bm{w})=\alpha$.
Equivalently, the line $t\mapsto t\cdot u(\bm{w})-v(\bm{w})$ is a supporting line for $\tau$ at $q$ with slope $\alpha$.
Therefore, the problem of determining the values of $\alpha$ for which $\tau$ is supported at $(q,\alpha)$ is precisely the problem of determining the slopes of supporting lines which appear from the minimization defining $\tau(q)$.

For the remainder of this section, we establish some elementary facts concerning the dual problems $\tau(q)$ and $f(\alpha)$.
We begin with the following basic fact about supports of the function $\tau$.
\begin{lemma}\label{l:supporting}
    For any $q\in\R$ and $\alpha\in\{\partial^-\tau(q),\partial^+\tau(q)\}$, $\tau$ is supported at $(q,\alpha)$.
\end{lemma}
\begin{proof}
    Let $(q_n)_{n=1}^\infty$ be a sequence converging to $q$ monotonically from the left.
    For each $n\in\N$, write $\tau(q_n)=q_n\cdot u(\bm{w}_n)-v(\bm{w}_n)$ for some $\bm{w}_n\in\Delta$.
    By compactness of $\Delta$, passing to a subsequence if necessary, we may assume that $\lim_{n\to\infty}\bm{w}_n=\bm{w}$.
    Next, observe that the line
    \begin{equation*}
        \ell_n(t)\coloneqq t\cdot u(\bm{w}_n)-v(\bm{w}_n)
    \end{equation*}
    is a supporting line for $\tau$ at $q_n$ and therefore $u(\bm{w}_n)\in\partial\tau(q_n)$.
    Moreover, since $\lim_{n\to\infty}\partial^-\tau(q_n)=\partial^-\tau(q)$, it follows that
    \begin{equation*}
        u(\bm{w})=\lim_{n\to\infty}u(\bm{w}_n)=\partial\tau^-(q).
    \end{equation*}
    But then by upper semicontinuity of $v$ and continuity of $\tau$,
    \begin{equation*}
        \tau(q)\leq q\cdot u(\bm{w})-v(\bm{w})\leq \lim_{n\to\infty}(q\cdot u(\bm{w}_n)-v(\bm{w}_n))=\lim_{n\to\infty}\tau(q_n)=\tau(q)
    \end{equation*}
    so that equality holds and $\tau$ is supported at $(q,\partial^-\tau(q))$.
    The same argument works for $\partial^+\tau(q)$, giving the result.
\end{proof}
\begin{remark}\label{r:bounded-derivs}
    Since $u$ is continuous and $\Delta$ is compact, it follows that the left and right derivatives are uniformly bounded away from $\pm\infty$.
\end{remark}
Using \cref{l:supporting}, we can now characterize concavity of the function $f$.
\begin{proposition}\label{p:conc-char}
    For any $\alpha\in\R$, $f(\alpha)\leq\tau^*(\alpha)$, and if $\tau$ is supported at $(q,\alpha)$ for some $q\in\R$, then $f(\alpha)=q\alpha-\tau(q)=\tau^*(\alpha)$.
    Moreover, the following are equivalent:
    \begin{enumerate}[nl,r]
        \item\label{im:supp} For all $q\in\R$ and $\alpha\in\partial\tau(q)$, $\tau$ is supported at $(q,\alpha)$.
        \item\label{im:conj} $f(\alpha)=\tau^*(\alpha)$ for all $\alpha\in\R$.
        \item\label{im:conc} $f$ is a concave function.
    \end{enumerate}
\end{proposition}
\begin{proof}
    First, let $\alpha\in\R$ be arbitrary.
    If $f(\alpha)=-\infty$, we are done; otherwise, since $\Delta$ is compact, $u$ is continuous, and $v$ is upper semi-continuous, there is some $\bm{w}\in\Delta$ so that $u(\bm{w})=\alpha$ and $v(\bm{w})=f(\alpha)$.
    Then for any $q\in\R$,
    \begin{equation*}
        \tau(q)\leq q\cdot u(\bm{w})-v(\bm{w})=q\alpha-f(\alpha).
    \end{equation*}
    But $q$ was arbitrary, so $f(\alpha)\leq\tau^*(\alpha)$.

    Next, if $\tau$ is supported at $(q,\alpha)$ for $q\in\R$, let $\bm{w}$ be such that $u(\bm{w})=\alpha$ and $\tau(q)=q\alpha-v(\bm{w})$.
    But then
    \begin{equation*}
        \tau^*(\alpha)\geq f(\alpha)\geq v(\bm{w})=q\alpha-\tau(q)\geq\tau^*(\alpha).
    \end{equation*}
    as claimed.

    Now, that \cref{im:supp} implies \cref{im:conj} was proven above, and \cref{im:conj} immediately implies \cref{im:conc}.
    It remains to verify that \cref{im:conc} implies \cref{im:supp}.
    Let $q\in\R$.
    Then if $\alpha\in\partial\tau(q)$, write $\alpha=\lambda \partial^+\tau(q)+(1-\lambda)\partial^-\tau(q)$ for some $\lambda\in[0,1]$.
    Then by concavity of $f$ and \cref{l:supporting},
    \begin{align*}
        \tau^*(\alpha)&\geq f(\alpha)\\
                      &\geq\lambda f\bigl(\partial^+\tau(q)\bigr)+(1-\lambda)f\bigl(\partial^-\tau(q)\bigr)\\
                      &=\lambda (q\partial^+\tau(q)-\tau(q))+(1-\lambda)(q\partial^-\tau(q)-\tau(q))\\
                      &=q\alpha-\tau(q)\\
                      &\geq\tau^*(\alpha).
    \end{align*}
    Thus all the inequalities are in fact equalities.
    In particular, taking $\bm{w}$ so that $u(\bm{w})=\alpha$ and $f(\alpha)=v(\bm{w})$, substituting this into the previous equation implies that $v(\bm{w})=q\cdot u(\bm{w})-\tau(q)$, as required.
\end{proof}
To round off our results concerning the concave conjugate relationship, let's describe what happens as $q$ diverges to $\pm\infty$.
We first extend the derivatives to $\pm\infty$ by defining
\begin{equation*}
    \partial^-\tau(\infty)=\lim_{q\to \infty}\frac{\tau(q)}{q}\qquad\text{and}\qquad\partial^+\tau(-\infty)=\lim_{q\to -\infty}\frac{\tau(q)}{q}.
\end{equation*}
That these limits exist and take finite values follows by monotonicity of the partial derivatives combined with the observation in \cref{r:bounded-derivs}.
Thus $\tau$ has affine asymptotes at $\pm\infty$ with slopes $\partial^-\tau(\infty)$ and $\partial^+\tau(-\infty)$.
We then say that $\tau$ is supported at $(\infty,\partial^-\tau(\infty))$ if there is a $\bm{w}$ so that the line $q\mapsto q\cdot u(\bm{w})-v(\bm{w})$ is the affine asymptote of $\tau$ at $\infty$.
The same definition holds at $-\infty$ as well.
Then, the same proof as \cref{l:supporting} (along a sequence $(q_n)_{n=1}^\infty$ diverging to $\pm\infty$) combined with a modified \cref{p:conc-char} yields the following result.
\begin{proposition}
    The function $\tau$ is supported at $(\infty,\partial^+\tau(-\infty))$ and $(-\infty,\partial^-\tau(\infty))$.
    In particular, $\tau^*(\partial^+\tau(-\infty))=f(\partial^+\tau(-\infty))$ and $\tau^*(\partial^-\tau(\infty))=f(\partial^-\tau(\infty))$.
\end{proposition}
We conclude this section with two explicit situations in which we can establish the concave conjugate relationship.
The first situation, even without any knowledge of the underlying optimization, occurs when $\tau$ is differentiable.
\begin{corollary}\label{c:diff-equiv}
    If $q\in\R$ is such that $\alpha=\tau'(q)$ exists, then $f(\alpha)=\tau^*(\alpha)$.
    In particular, $\tau(q)=f^*(q)$ for all $q\in\R$.
\end{corollary}
\begin{proof}
    If $q\in\R$ and $\tau'(q)=\alpha$ exists, by \cref{l:supporting}, $\tau$ is supported at $q$.
    Thus by \cref{p:conc-char}, $f(\alpha)=\tau^*(\alpha)$.
    In particular, $\tau(q)=f^*(q)$ for all $q\in\R$ for which $\tau'(q)$ exists.
    But $f^*$ and $\tau$ are both continuous functions and $\tau'(q)$ exists for a dense set of $q\in\R$, so in fact $\tau(q)=f^*(q)$ for all $q\in\R$.
\end{proof}
As our second (and final) application, we can also use information about the structure of the set on which the optimization is attained to abstractly establish the concave conjugate relationship.
For each $q\in\R$, denote the set of minimizers by
\begin{equation*}
    M(q) = \{\bm{w}\in\Delta:q\cdot u(\bm{w})-v(\bm{w})=\tau(q)\}.
\end{equation*}
Since $u$ is continuous and $v$ is upper semicontinuous, $M(q)$ is a compact set for all $q$.
Moreover, we obtain the following result.
\begin{corollary}\label{c:path-conn}
    Suppose $q\in\R$ and $M(q)$ is connected.
    Then $f(\alpha)=\tau^*(\alpha)$ for all $\alpha\in\partial\tau(q)$.
    Moreover, if $M(q)$ is a singleton, then $\tau$ is differentiable at $q$.
\end{corollary}
\begin{proof}
    By \cref{l:supporting}, $\tau$ is supported at $(q,\partial^-\tau(q))$ and $(q,\partial^+\tau(q))$.
    Thus there are $\bm{w}_-,\bm{w}_+\in M(q)$ so that
    \begin{equation*}
        u(\bm{w}_{-})=\partial^{-}\tau(q)
        \qquad\text{and}\qquad
        u(\bm{w}_{+})=\partial^{+}\tau(q).
    \end{equation*}
    Since $M(q)$ is connected and $u$ is continuous, $u(M(q))\subset\R$ is an interval containing $\partial^-\tau(q)$ and $\partial^+\tau(q)$, and therefore $\partial\tau(q)\subset u(M(q))$.
    In particular, for any $\alpha\in\partial\tau(q)$, there is a $\bm{w}\in M(q)$ such that $u(\bm{w})=\alpha$.
    Thus $\tau$ is supported at $(q,\alpha)$ and therefore $f(\alpha)=\tau^*(\alpha)$ by \cref{p:conc-char}.

    If moreover $M(q)$ is a singleton, then we must have $\bm{w}_-=\bm{w}_+$, forcing $\partial^-\tau(q)=\partial^+\tau(q)$ so that $\tau$ is differentiable at $q$.
\end{proof}

\subsection{Multifractal spectrum and general upper bound}\label{ss:Lq-bound}
For $\alpha\in\R$, let
\begin{equation*}
    F_\mu(\alpha)=\{x\in\supp\mu:\dim_{\loc}(\mu,x)=\alpha\}.
\end{equation*}
Then the \defn{multifractal spectrum} is the function $f_\mu\colon\R\to\{-\infty\}\cup[0,d]$ given by
\begin{equation*}
    f_\mu(\alpha)=\dimH F_\mu(\alpha)
\end{equation*}
using the convention that $\dimH\varnothing=-\infty$.
We say that the measure $\mu$ satisfies the \defn{multifractal formalism} if
\begin{equation*}
    f_\mu(\alpha)=\tau_\mu^*(\alpha)
\end{equation*}
for all $\alpha\in\R$.

We first establish the general upper bound for the multifractal spectrum from the $L^q$-spectrum.
Such a bound has been known for a long time; see, for example, \cite[Theorem~4.1]{zbl:0929.28007}.
However, many of the proofs in the literature do not precisely address the cases when $\alpha$ corresponds to a slope of the asymptote of $\tau_\mu$.
We clarify this case explicitly in the following result.
\begin{proposition}\label{p:lq-concave-bound}
    Let $\mu$ be a compactly supported Borel probability measure.
    Then
    \begin{equation*}
        f_\mu(\alpha)\leq\tau_\mu^*(\alpha)
    \end{equation*}
    for all $\alpha\in\R$.
\end{proposition}
\begin{proof}
    For $\alpha\in\R$, $r\in(0,1)$ and $\varepsilon>0$, let
    \begin{equation*}
        \mathcal{M}_{r,\varepsilon}(\alpha)=\Bigl\{x\in\supp\mu:r^{\alpha+\varepsilon}\leq\mu\bigl(B(x,r)\bigr)\leq r^{\alpha-\varepsilon}\Bigr\}.
    \end{equation*}
    Our strategy is to control the size of packings with centres in $\mathcal{M}_{r,\varepsilon}(\alpha)$ using the $L^q$-spectrum of $\mu$, to which we can then apply the Hausdorff dimension version of the Vitali covering theorem (see \cite[Theorem 1.10]{zbl:0587.28004}) to obtain bounds on $f_\mu(\alpha)$.

    More precisely, we first prove the following claim.
    \begin{claim}\label{c:disj-ct}
        For all $\alpha\in\R$ and $\zeta>0$, there are $\varepsilon > 0$ and $r_0 > 0$ (depending only on $\zeta$ and $\alpha$) such that for any $0<r \leq r_0$ and any disjoint family of balls $\{B(x_i,r)\}_{i=1}^N$ with $x_i\in \mathcal{M}_{r,\varepsilon}$,
        \begin{equation}\label{e:M-cover-bound}
            \frac{\log N}{\log(1/r)}\leq \tau_\mu^*(\alpha)+\zeta.
        \end{equation}
    \end{claim}
    Fix $\varepsilon > 0$ and $r > 0$ (we will constrain these values below), and fix an arbitrary disjoint family of balls $\{B(x_i,r)\}_{i=1}^N$ where $x_i\in \mathcal{M}_{r,\varepsilon}$.
    Suppose $q\geq 0$ is arbitrary.
    Then (recalling the definition of $G_\mu$ from the start of \cref{ss:lq-spectrum})
    \begin{equation}\label{e:sn-upper}
        G_\mu(r,q)\geq\sum_{i=1}^N\mu\bigl(B(x_i,r)\bigr)^q\geq N\cdot r^{q(\alpha+\varepsilon)}.
    \end{equation}
    Moreover, since $\tau_\mu(q)=\liminf_{r\to 0}(\log G_\mu(r,q))/(\log r)$, there is some $r_\varepsilon\in(0,1)$ so that for all $r\leq r_\varepsilon$, $G_\mu(r,q)\leq r^{\tau_\mu(q)-\varepsilon}$.
    Combining this with \cref{e:sn-upper},
    \begin{equation}\label{e:M-q-bound}
        \frac{\log N}{\log(1/r)}\leq q\alpha-\tau_\mu(q)+(|q|+1)\varepsilon
    \end{equation}
    for all $r\leq r_\varepsilon$.
    The same argument for $q<0$ also yields \cref{e:M-q-bound} using the bound $\mu(B(x_i,r))\leq r^{\alpha-\varepsilon}$.

    Now we verify \cref{e:M-cover-bound}.
    To do this, we must handle a few different cases depending on the choice of $\alpha$.
    Recall that $\tau_\mu(q)$ is concave with asymptotes at $\pm\infty$ given respectively by
    \begin{equation*}
        q\mapsto q\cdot a_{-\infty}-\tau_\mu^*(a_{-\infty})\qquad\text{and}\qquad q\mapsto q\cdot a_{\infty}-\tau_\mu^*(a_\infty).
    \end{equation*}
    If $a_{\infty}<\alpha<a_{-\infty}$, there is a $q\in\R$ so that $q\in\partial\tau_\mu^*(\alpha)$ and therefore $q\alpha-\tau_\mu(q)=\tau^*(\alpha)$.
    Taking $\varepsilon$ sufficiently small so that $(|q|+1)\varepsilon\leq\zeta$ then yields the claim.
    Next, if $\alpha=a_\infty$ (the case $\alpha=a_{-\infty}$ is analogous), since
    \begin{equation*}
        \lim_{q\to\infty}(q\alpha-\tau(q))=\tau_\mu^*(\alpha),
    \end{equation*}
    we first get $q$ so that $q\alpha-\tau(q)\leq\tau_\mu^*(\alpha)+\zeta/2$, and then get $\varepsilon>0$ so that $(|q|+1)\varepsilon\leq\zeta/2$, as claimed.
    Finally, if $\alpha<-a_{\infty}$ or $\alpha>a_\infty$, taking an infimum over all $q\in\R$ yields $N=0$ for $\varepsilon$ sufficiently small depending on $\alpha$.
    This completes the proof of \cref{c:disj-ct}.

    \vspace{0.25cm}

    Now for each $x\in F_\mu(\alpha)$, we can find some $N_x\in\N$ so that for all $n\geq N_x$, $\mu\bigl(B(x,2^{-n})\bigr)\geq 2^{-n(\alpha+\varepsilon)}$.
    In particular, for any $m\in\N$,
    \begin{equation*}
        \mathcal{G}_{m, \varepsilon}\coloneqq\bigcup_{n=m}^\infty\mathcal{M}_{2^{-n},\varepsilon}(\alpha)
    \end{equation*}
    is a Vitali cover for $F_\mu(\alpha)$.

    To conclude the proof, fix $\alpha\in\R$, let $\xi>0$ be arbitrary and get $m\in\N$ so that \cref{e:M-cover-bound} holds for all $r\leq 2^{-m}$.
    Suppose $\{B(x_j,2^{-n_j})\}_{j=1}^\infty$ is a disjoint subcollection of $\mathcal{G}_{m,\varepsilon}$: then with $s=\tau^*(\alpha)+2\xi$,
    \begin{equation*}
        \sum_{j=1}^\infty 2^{-n_j s}\leq\sum_{n=m}^\infty 2^{-n s}2^{-n(\tau^*(\alpha)+\xi)}=\sum_{n=m}^\infty (2^{-\xi})^{n}<\infty
    \end{equation*}
    by \cref{c:disj-ct}.
    Thus by the Vitali covering theorem for Hausdorff measure, there is a cover $\{E_i\}_{i=1}^\infty$ for $F_\mu(\alpha)$ such that
    \begin{equation*}
        \mathcal{H}^s(F_\mu(\alpha))\leq \sum_{i=1}^\infty(\diam E_i)^s<\infty
    \end{equation*}
    and thus $f_\mu(\alpha)\leq\tau^*(\alpha)+2\xi$.
    Since $\xi>0$ was arbitrary, the result follows.
\end{proof}
\begin{remark}
    Why is this bound not sharp in general?
    In the proof, we are controlling the number of balls at resolution $r$ which have measure larger (or smaller) than the expected quantity.
    However, the set of points with local dimension $\alpha$ might always be hard to cover at some fixed scale $r$, even though some parts of the set are easier to cover at some scales than others.
    This is similar to the situation in which the Hausdorff dimension and lower box dimension differ---which is precisely the case here when $q=0$.

    In fact, one might think of the concave conjugate $\tau_\mu^*(\alpha)$ as a \emph{geometric large deviations bound} for the multifractal spectrum.
    Instead of bounding the set of points with non-typical scaling, the $L^q$-spectrum bounds the asymptotic size of the set of points with non-typical scaling at a fixed scale $r$, as $r$ converges to zero.
    This analogy is particularly relevant for a measure with dynamical origin (this is highlighted, for instance, in the proof of \cref{t:lq-form}).
\end{remark}
\subsection{Multifractal formalism for self-similar measures}
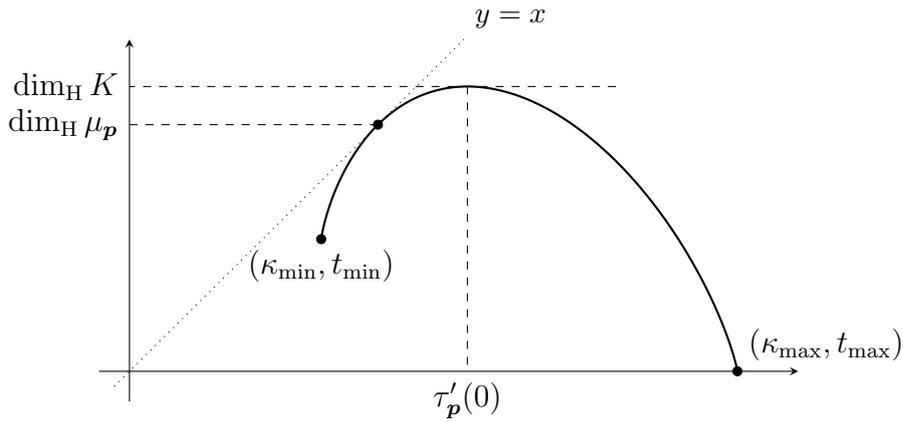
\begin{figure}[t]
    \centering
    \input{figures/multifractal}
    \caption{Plot of the multifractal spectrum $f_{\bm{p}}(\alpha)$.}
    \label{f:multifractal}
\end{figure}

We can now establish the multifractal formalism for self-similar measures satisfying the SSC by combining the various tools that we introduced in the previous sections.
Again, for short, given a self-similar measure $\mu_{\bm{p}}$ we write $f_{\mu_{\bm{p}}}=f_{\bm{p}}$.
A depiction of the multifractal spectrum $f_{\bm{p}}$ is given in \cref{f:multifractal}, with the same parameters from \cref{f:lq}.
The formula for the multifractal spectrum $f_{\bm{p}}$ is originally due to Cawley \& Mauldin \cite{zbl:0763.58018}, and the relationship with the $L^q$-spectrum established by Arbeiter \& Patzschke \cite{zbl:0873.28003}.
Also see the earlier closely related work due to Rand \cite{zbl:0664.58022}.
\begin{theorem}[\cite{zbl:0763.58018,zbl:0873.28003}]
    Let $\{S_i\}_{i\in\mathcal{I}}$ be an IFS satisfying the SSC and let $\bm{p}\in\mathcal{P}$.
    Then $\mu_{\bm{p}}$ satisfies the multifractal formalism.
\end{theorem}
\begin{proof}
    A direct application of \cref{p:rel-dim} gives that $f_{\bm{p}}\geq f$ where
    \begin{equation*}
        f(\alpha)\coloneqq\sup_{\bm{w}\in\mathcal{P}}\left\{\frac{H(\bm{w})}{\chi(\bm{w})}:\frac{H(\bm{w},\bm{p})}{\chi(\bm{w})}=\alpha\right\}.
    \end{equation*}
    Moreover, by \cref{t:lq-form} and \cref{c:diff-equiv} combined with \cref{p:tau-diff}, $f=\tau_{\bm{p}}^*$.
    Finally, $f_{\bm{p}}\leq\tau_{\bm{p}}^*$ holds in general by \cref{p:lq-concave-bound}, giving the desired result.
\end{proof}
\begin{remark}
    Alternatively, one can use \cref{p:var-upper} in the following section, in place of the general upper bound from \cref{p:lq-concave-bound}.
    In particular, if one only wishes to obtain the dimensional formula for $f_{\bm{p}}(\alpha)$, this provides a method to entirely bypass the results concerning $L^q$-spectra.
\end{remark}
\begin{remark}
    In our proof of the multifractal formalism, we obtain a formula for $f_{\bm{p}}$ which also includes the endpoints.
    Many proofs of the multifractal formalism either ignore the behaviour at the endpoint entirely, or handle it by explicit arguments which are qualitatively different than the main proof of the multifractal formalism.
    In contrast, this optimization technique gives the result at the endpoint for free: the endpoint behaviour is simply attained as the limit (in the compact space $\mathcal{P}$) of the optimizing vector for the $L^q$-spectrum.
    The geometry of this parametrization is discussed in more detail in \cref{s:geom}.
\end{remark}
\subsection{Alternative proof of the variational principle}\label{ss:alt}
In this section, we prove a general covering theorem for the Hausdorff dimension of sets.
This bound is similar to the general upper bound for the Hausdorff dimension using measures (see, for example, \cite[Proposition~10.1]{zbl:0869.28003}), though we will state and prove an extension which permits the choice of measure to depend on both the point in the set and the scale.
This additional flexibility will allow us to provide an alternative proof of the variational formula for the level sets of local dimensions.
\begin{definition}
    Let $\Delta$ be a family of Borel probability measures on $\R^d$.
    We say that $\Delta$ has \defn{uniform densities} if for every $\varepsilon>0$, there exists $C>0$, $\eta\in(0,1)$, and a family $\mathcal{E}$ of Borel measures with $\sum_{\nu\in\mathcal{E}}\nu(\R^d)<\infty$ so that for all $\mu\in\Delta$, $r\in(0,\eta)$, and $x\in\R^d$, there is a $\nu=\nu_{\mu,x,r}\in\mathcal{E}$ so that
    \begin{equation*}
        r^\varepsilon\mu\bigl(B(x,r)\bigr)\leq\nu\bigl(B(x,Cr)\bigr).
    \end{equation*}
\end{definition}
Intuitively, having uniform densities combines a uniform semi-continuity condition on the local dimensions of measures in $\Delta$ with a pre-compactness condition on $\Delta$.
We emphasize that we do not require the measures in $\mathcal{E}$ to belong to $\Delta$.

A trivial example is any finite set of Borel probability measures.
Less trivially, and more usefully for us, the set of Bernoulli measures associated with a self-similar IFS (not necessarily satisfying the SSC) also has uniform densities.
\begin{proposition}\label{p:P-compact}
    Let $\{S_i\}_{i\in\mathcal{I}}$ be a self-similar IFS with attractor $K$.
    Then $\{\mu_{\bm{p}}:\bm{p}\in\mathcal{P}\}$ has uniform densities.
\end{proposition}
\begin{proof}
    Let $\xi = \max\{r_i:i\in\mathcal{I}\}$ and let $\varepsilon>0$ be arbitrary.
    For each $\bm{z}\in\mathcal{P}$ with $\bm{z} > 0$, let
    \begin{equation*}
        U(\bm{z}) = \{\bm{w}\in\mathcal{P}: \xi^\varepsilon w_i < z_i\text{ for all }i\in\mathcal{I}\}.
    \end{equation*}
    Since $\bm{z} > 0$, it is easy to check that $U(\bm{z})$ is an open set.
    Let us check that $\{U(\bm{z}):\bm{z} > 0\}$ is in fact a cover for $\mathcal{P}$.
    Let $\bm{w}\in\mathcal{P}$ be arbitrary and let $\mathcal{I}_0 = \{i\in\mathcal{I}: w_i = 0\}$.
    For all $t > 0$ sufficiently small, define the probability vector $\bm{z}(t)$
    \begin{equation*}
        z_i(t) = \begin{cases}
            t &: i \in \mathcal{I}_0,\\
            w_i - t\frac{\#\mathcal{I}_0}{\#\mathcal{I} - \#\mathcal{I}_0} &: i \notin\mathcal{I}_0.
        \end{cases}
    \end{equation*}
    Moreover, it is easy to check for all $t$ sufficiently small depending on $\bm{w}$ and $\xi^{\varepsilon}$ that $\bm{w} \in U(\bm{z}(t))$.
    Since $\mathcal{P}$ is compact, there exists a finite subcover
    \begin{equation*}
        \{U(\bm{z}_1),\ldots,U(\bm{z}_n)\}
    \end{equation*}
    of $\mathcal{P}$.

    Let us check that $\mathcal{E} = \{\mu_{\bm{z}_1},\ldots,\mu_{\bm{z}_n}\}$ satisfies the required conditions in the definition of uniform densities.
    For each $r\in(0,1)$ and $x\in\R^d$, let
    \begin{equation*}
        \Lambda(x,r) = \bigl\{\mtt{i}\in\mathcal{I}^m: \xi^{m+1} < r \leq \xi^{m}, S_{\mtt{i}}(K) \cap B(x,r)\neq\varnothing\bigr\}.
    \end{equation*}
    If $\mtt{i}\in\Lambda(x,r)$, then $\diam S_{\mtt{i}}(K) \leq r \xi^{-1}\diam K$ so
    \begin{equation}\label{e:C-choice}
        S_{\mtt{i}}(K) \subset B(x, Cr)\quad\text{where}\quad C \coloneqq 1 + \xi^{-1}\diam K.
    \end{equation}
    Let $\bm{w}\in\mathcal{P}$ be arbitrary and let $\bm{z}\in\{\bm{z}_1,\ldots,\bm{z}_n\}$ be such that $\bm{w}\in U(\bm{z})$.
    Then
    \begin{align*}
        r^{\varepsilon} \mu_{\bm{w}}\bigl(B(x,r)\bigr) &\leq \sum_{\mtt{i}\in\Lambda(x,r)} r^{\varepsilon}w_{\mtt{i}}\\
                                              &\leq \sum_{\mtt{i}\in\Lambda(x,r)}\prod_{j=1}^m(\xi^\varepsilon w_{\mtt{i}_j})&&\text{since $r \leq \xi^m$}\\
                                              &\leq \sum_{\mtt{i}\in\Lambda(x,r)} z_{\mtt{i}}&&\text{since $\bm{w}\in U(\bm{z})$}\\
                                              & \leq \mu_{\bm{z}}\bigl(B(x, Cr)\bigr)&&\text{by \cref{e:C-choice}}
    \end{align*}
    as required.
\end{proof}
In the following result, the case where $\Delta$ is a singleton is standard: the additional flexibility where the measure $\mu$ can be chosen to depend on both the point and the scale is useful in applications.
\begin{proposition}\label{p:H-bound}
    Let $K\subset\R^d$ be an arbitrary set and let $\Delta$ be a family of measures with uniform densities.
    Suppose $t$ is such that for all $x\in K$,
    \begin{equation*}
        \liminf_{r\to 0}\inf_{\mu\in\Delta}\frac{\log\mu\bigl(B(x,r)\bigr)}{\log r}\leq t.
    \end{equation*}
    Then $\dimH K\leq t$.
\end{proposition}
\begin{proof}
    By countable stability of Hausdorff dimension, it suffices to consider bounded $K$.
    Fix $\varepsilon>0$.
    Since $\Delta$ has uniform densities, there are constants $C>0$, $\eta_0\in(0,1)$ and a family $\mathcal{E}$ of finite Borel measures on $\R^d$ with $\sum_{\nu\in\mathcal{E}}\nu(\R^d)<\infty$ so that for any $\mu\in\Delta$, $x\in K$, and $r\in(0,\eta_0)$, there is a $\nu=\nu_{\mu,x,r}\in\mathcal{E}$ so that
    \begin{equation}\label{e:discretize}
        r^\varepsilon\mu\bigl(B(x,r)\bigr)\leq \nu\bigl(B(x, C r)\bigr).
    \end{equation}
    For each $\eta\in(0,\eta_0)$ and $x\in K$, choose a measure $\mu_x\in\Delta$ and a radius $r_x\in(0,\eta)$ such that
    \begin{equation}\label{e:mux-choice}
        r_x^{t+\varepsilon}\leq\mu_x\bigl(B(x,r_x)\bigr).
    \end{equation}

    Let $\mathcal{B}=\{B(x, C r_x)\}_{x\in K}$.
    By the Besicovitch covering theorem (see \cite[Theorem~2.7]{zbl:0819.28004}), there is a constant $c_d$ depending only on the ambient dimension $d$ and families of balls $\mathcal{B}_i\subset\mathcal{B}$ for $i=1,\ldots,c_d$ such that
    \begin{equation*}
        K\subset\bigcup_{i=1}^{c_d}\bigcup_{B\in\mathcal{B}_i}B
    \end{equation*}
    and moreover the balls in each $\mathcal{B}_i$ are disjoint.
    Thus for each $i=1,\ldots,c_d$ and $B(x,Cr)\in\mathcal{B}_i$, by \cref{e:mux-choice} combined with \cref{e:discretize}, there is a $\nu\in\mathcal{E}$ so that
    \begin{equation*}
        r^{t+2\varepsilon}\leq\nu\bigl(B(x,C r)\bigr).
    \end{equation*}
    But the balls in each $\mathcal{B}_i$ are pairwise disjoint, so
    \begin{align*}
        \sum_{i=1}^{c_d}\sum_{B\in\mathcal{B}_i}(\diam B)^{t+2\varepsilon}
        &\leq \sum_{i=1}^{c_d}\sum_{B(x,Cr)\in\mathcal{B}_i}(2C r)^{t+2\varepsilon}\\
        &\leq(2C)^{t + 2\varepsilon}\sum_{i=1}^{c_d}\sum_{\nu\in\mathcal{E}}\nu\left(\bigcup_{B\in\mathcal{B}_i}B\right)\\
        &\leq c_d(2C)^{t+2\varepsilon}\cdot\left(\sum_{\nu\in\mathcal{E}}\nu(\R^d)\right)<\infty
    \end{align*}
    which, since $\eta>0$ and $\varepsilon>0$ were arbitrary, gives that $\dimH K\leq t$.
\end{proof}
We now obtain the following variational upper bound.
In fact, our result is a bit stronger since we can actually bound any point such that the local dimension $\alpha$ is attained along a subsequence.
More precisely, given a Borel probability measure $\mu$ and $x\in \supp\mu$, we write
\begin{equation*}
    \mathcal{D}(\mu,x)=\left\{\alpha:\exists(r_n)_{n=1}^\infty\searrow 0\text{ s.t.\ }\lim_{n\to\infty}\frac{\log\mu(B(x,r_n))}{\log r_n}=\alpha\right\}.
\end{equation*}
Of course, $\dim_{\loc}(\mu,x)=\alpha$ if and only if $\mathcal{D}(\mu,x)=\{\alpha\}$.
The following formula was first established in \cite{zbl:1040.28014} via a direct covering argument.
\begin{proposition}[\cite{zbl:1040.28014}]\label{p:var-upper}
    Let $\{S_i\}_{i\in\mathcal{I}}$ be a self-similar IFS satisfying the SSC, let $\bm{p}\in\mathcal{P}$, and let $\alpha\in\R$.
    Then
    \begin{equation}\label{e:f-pressure-form}
        \begin{aligned}
            f_{\bm{p}}(\alpha)&=\dimH\{x\in K:\alpha\in\mathcal{D}(\mu_{\bm{p}},x)\}\\
                              &=\sup_{\bm{w}\in\mathcal{P}}\left\{\frac{H(\bm{w})}{\chi(\bm{w})}:\frac{H(\bm{w},\bm{p})}{\chi(\bm{w})}=\alpha\right\}.
        \end{aligned}
    \end{equation}
\end{proposition}
\begin{proof}
    Recalling \cref{p:rel-dim}, it suffices to prove that
    \begin{equation*}
        \dimH\{x\in K:\alpha\in\mathcal{D}(\mu_{\bm{p}},x)\}\leq\sup_{\bm{w}\in\mathcal{P}}\left\{\frac{H(\bm{w})}{\chi(\bm{w})}:\frac{H(\bm{w},\bm{p})}{\chi(\bm{w})}=\alpha\right\}.
    \end{equation*}
    To this end, fix $\gamma\in\Omega$ and $\alpha\in\mathcal{D}(\mu_{\bm{p}},x)$ where $x=\pi(\gamma)$.
    Write $\bm{q}_n=\bm{\xi}(\gamma\npre{n})$ and $R_n=D r_{\gamma\npre{n}}$, where $D>\max\{1,\diam K\}$ is fixed.
    By (the proof of) \cref{l:symb-formula}, there is a subsequence $(n_k)_{k=1}^\infty$ so that $\lim_{k\to\infty}\bm{q}_{n_k}=\bm{w}$ and
    \begin{equation*}
        \lim_{k\to\infty}\frac{\log \mu_{\bm{p}}\bigl(B(x,r_{\gamma\npre{n_k}})\bigr)}{\log r_{\gamma\npre{n_k}}}=\alpha.
    \end{equation*}
    In particular,
    \begin{equation*}
        \frac{H(\bm{w},\bm{p})}{\chi(\bm{w})}=\alpha.
    \end{equation*}
    Since $S_{\gamma\npre{n}}(K)\subset B(x,R_n)$, we have for every $n$
    \begin{equation}\label{e:empirical-cylinder}
        \log \mu_{\bm{q}_n}\bigl(B(x,R_n)\bigr)\geq -nH(\bm{q}_n).
    \end{equation}
    Consequently, continuity of entropy on $\mathcal{P}$ gives
    \begin{equation*}
        \liminf_{r\to 0}\inf_{\bm{v}\in\mathcal{P}}\frac{\log\mu_{\bm{v}}\bigl(B(x,r)\bigr)}{\log r}\leq\lim_{k\to\infty}\frac{n_k H(\bm{q}_{n_k})}{-\log R_{n_k}}=\frac{H(\bm{w})}{\chi(\bm{w})}.
    \end{equation*}
    By \cref{p:P-compact}, $\mathcal{P}$ has uniform densities.
    Thus the Hausdorff dimension bound follows by \cref{p:H-bound}.
\end{proof}
\begin{remark}
    As discussed in the introduction, both the technique of uniform densities, as well as the method of types, are useful tools for reducing a general covering argument to a covering argument for ``typical'' points of a well-chosen measure.
    Heuristically, one might interpret the method of types argument as the ``box-dimension'' variant of this reduction, whereas the uniform densities argument is the ``Hausdorff'' variant of this reduction.
\end{remark}
\begin{remark}
    A similar argument as presented in this section will allow us to instead bound the \emph{packing} dimension of the set $\{x\in K:\dim_{\loc}(\mu_{\bm{p}},x)=\alpha\}$.
    In order to bound the packing dimension, in \cref{p:H-bound}, instead of providing a cover at infinitely many scales, one must provide a cover at \emph{all} scales.
    More precisely, if one replaces the limit infimum with a limit supremum in \cref{p:H-bound}, one instead obtains a bound for the packing dimension of the underlying set.
    Then a similar argument as given in \cref{p:var-upper} shows that
    \begin{equation*}
        \dimP\{x\in K:\alpha\in\mathcal{D}(\mu_{\bm{p}},x)\}=\alpha\leq\sup_{\bm{w}\in\mathcal{P}}\left\{\frac{H(\bm{w})}{\chi(\bm{w})}:\frac{H(\bm{w},\bm{p})}{\chi(\bm{w})}=\alpha\right\}.
    \end{equation*}
    The key difference is that the subsequence $(n_k)_{k=1}^\infty$ must now be chosen to realize the limit supremum, which necessitates the assumption that the local dimension exists.
    In fact, it is no longer true that the packing dimension of the set $\{x\in K:\alpha\in\mathcal{D}(\mu_{\bm{p}},x)\}$ satisfies the same bound (see, for instance, \cite{zbl:1141.28006}).
\end{remark}
\section{Concluding thoughts: The geometry of optimization}\label{s:geom}
\subsection{Recalling duality}
Recall the setting from \cref{ss:optim}: $\Delta$ is a compact set and we are given a continuous function $u\colon\Delta\to\R$ and an upper semicontinuous function $v\colon\Delta\to\R$.
In fact, we now have a bit more context since in the multifractal analysis of self-similar measures, $\Delta=\mathcal{P}$ is the set of probability vectors, the function $v$ is the dimension of the corresponding pushforward measure
\begin{equation*}
    v(\bm{w})=\dimH\mu_{\bm{w}}=\frac{H(\bm{w})}{\chi(\bm{w})}
\end{equation*}
and the function $u$ gives the $\mu_{\bm{w}}$-typical local dimensions of $\mu_{\bm{p}}$
\begin{equation*}
    u(\bm{w})=\frac{H(\bm{w},\bm{p})}{\chi(\bm{w})}.
\end{equation*}
We then considered the constrained optimization problem
\begin{equation*}
    f(\alpha)=\max_{\bm{w}\in\Delta}\bigl\{v(\bm{w}):u(\bm{w})=\alpha\bigr\}
\end{equation*}
with unconstrained dual
\begin{equation*}
    \tau(q)=\min_{\bm{w}\in\Delta}\bigl\{q\cdot u(\bm{w})-v(\bm{w})\bigr\}.
\end{equation*}
In our particular situation $f(\alpha)=f_{\bm{p}}(\alpha)$ is precisely the multifractal spectrum, and $\tau(q)=\tau_{\bm{p}}(q)$ is the $L^q$-spectrum.

In \cref{p:tau-diff}, we proved that $\sum_{i\in\mathcal{I}}p_i^q r_i^{-\tau_{\bm{p}}(q)}=1$.
In particular, by the analytic implicit function theorem, this gives that $\tau_{\bm{p}}$ is differentiable; this was used in conjunction with \cref{c:diff-equiv} to prove the concave conjugate relationship.
In fact, we proved more: we showed that
\begin{equation*}
    \tau_{\bm{p}}(q)=q\cdot u(\bm{z}(q))-v(\bm{z}(q))
\end{equation*}
where $\bm{z}(q)=\bigl(p_i^q r_i^{-\tau_{\bm{p}}(q)}\bigr)_{i\in\mathcal{I}}$ defines a continuous path in $\mathcal{P}$.
Moreover, the minimization giving $\tau_{\bm{p}}(q)$ is attained uniquely at the vector $\bm{z}(q)$.

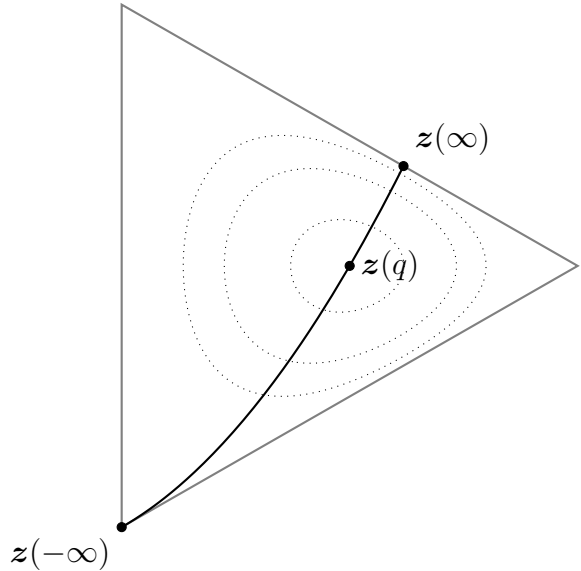
\begin{figure}[t]
    \centering
    \input{figures/simplex}
    \caption{Parametrization of the optimization vector in the simplex $\mathcal{P}\subset\R^3$ and contour lines of $\bm{w}\mapsto q\cdot u(\bm{w})-v(\bm{w})$.}
    \label{f:simplex}
\end{figure}

We note that one can also implicitly verify that the minimization is attained uniquely at a vector $\bm{z}(q)$.
Fix some $q\in\R$ and for each $t\in\R$, consider the lower level set
\begin{equation*}
    E(t,q)=\left\{\bm{w}\in\mathcal{P}:\frac{q\cdot H(\bm{w},\bm{p})-H(\bm{w})}{\chi(\bm{w})}\leq t\right\}.
\end{equation*}
Rearranging the condition on $E(t,q)$, we equivalently see that
\begin{equation*}
    E(t,q)=\bigl\{\bm{w}\in\mathcal{P}:H(\bm{w})-q\cdot H(\bm{w},\bm{p}) +t\cdot\chi(\bm{w})\geq 0\bigr\}.
\end{equation*}
But $H$ is a strictly convex function and $\bm{w}\mapsto t\cdot\chi(\bm{w})-q\cdot H(\bm{w},\bm{p})$ is a linear function in $\bm{w}$, so $E(t,q)$ is a convex set.
Therefore, noting that
\begin{equation*}
    \tau_{\bm{p}}(q)=\inf\{t:E(t,q)\neq\varnothing\},
\end{equation*}
we must have that $E(\tau_{\bm{p}}(q),q)=\{\bm{z}(q)\}$.
Since all the relevant functions are continuous in $q$, we also conclude that $\bm{z}$ is continuous in $q$.
The simplex in $\R^3$, the curve $\bm{z}(q)$, and the level sets of the function $\bm{w}\mapsto q\cdot u(\bm{w})-v(\bm{w})$, are depicted in \cref{f:simplex}.
The parameters are the same as the parameters used in \cref{f:lq} and \cref{f:multifractal},

This technique is useful to determine uniqueness of the minimization, or continuous path optimization, even in situations where it may be difficult to determine an explicit formula for the minimization.
In fact, we just proved that the objective function $\phi(\bm{w})= q\cdot u(\bm{w})-v(\bm{w})$ is \emph{strictly quasiconvex}: that is, for all $\bm{w}_1,\bm{w}_2\in\mathcal{P}$ and $\lambda\in(0,1)$,
\begin{equation*}
    \phi\bigl(\lambda\bm{w}_1+(1-\lambda)\bm{w}_2\bigr)<\max\{\phi(\bm{w}_1),\phi(\bm{w}_2)\}.
\end{equation*}
We say that $\phi$ is \emph{quasiconvex} if the above inequality holds with a regular inequality in place of a strict inequality.
Quasiconvexity is a useful property since it follows directly from the definition that for continuous quasiconvex functions defined on a compact convex set, the minimum must be attained on a convex subset.
Moreover, strict quasiconvexity implies that the minimum must be attained at a singleton.
In particular, (strict) quasiconvexity of the objective function $\bm{w}\mapsto q\cdot u(\bm{w})-v(\bm{w})$ at some $q\in\R$ allows us to immediately conclude the validity of the variational formula (or even differentiability) at $q$, by using \cref{c:path-conn}.

While quasiconvex functions share many properties with convex functions, one important exception is that a sum of quasiconvex functions is \emph{not} necessarily quasiconvex.
For instance, this is why the higher-dimensional versions of the optimization problems discussed in this paper are poorly behaved in general, such as non-uniqueness of measures of maximal dimensions \cite{zbl:1228.37025}.

\subsection{Parametrizing the optimization}\label{ss:param}
It turns out that even abstractly, having a continuously parametrized family of minimizers yields differentiability of $\tau$ along with other geometric information.
Suppose there is a continuous function $\bm{z}\colon\R\to\Delta$ such that $\tau(q)= q\cdot u(\bm{z}(q))-b(\bm{z}(q))$ for all $q\in\R$.
Using compactness of $\Delta$, we extend $\bm{z}$ to take values at $\pm\infty$ by taking any limit along a subsequence of $q$ diverging to $\pm\infty$.
We establish the following properties.
\begin{proposition}
    Suppose $\bm{z}\colon\R\to\Delta$ satisfies the above conditions.
    Then the following hold:
    \begin{enumerate}[nl,r]
        \item\label{im:deriv} The function $\tau$ is differentiable with derivative $\tau'(q)=u(\bm{z}(q))$.
        \item\label{im:decr} The function $q\mapsto u(\bm{z}(q))$ is monotonically decreasing.
        \item\label{im:fa} If $\alpha\in[u(\bm{z}(\infty)),u(\bm{z}(-\infty))]$, then
            \begin{equation*}
                f(\alpha)=v(\bm{z}(q_\alpha))
            \end{equation*}
            where $q_\alpha\in\R$ is any value so that $u(\bm{z}(q_\alpha))=\alpha$.
        \item\label{im:outside-dom} If $\alpha\notin[u(\bm{z}(\infty)),u(\bm{z}(-\infty))]$, then $f(\alpha)=-\infty$.
    \end{enumerate}
\end{proposition}
\begin{proof}
    First, to see \cref{im:deriv}, for each $q\in\R$, since $\bm{z}(q)$ attains the minimum for $\tau(q)$, it follows that $u(\bm{z}(q))\in\partial\tau(q)$.
    Moreover, by concavity of $\tau(q)$, if $q_1<q_2$, then $\inf\partial\tau(q_1)>\sup\partial\tau(q_2)$.
    Since $\bm{z}$ is continuous, this forces $\partial\tau(q)=\{u(\bm{z}(q))\}$, so that $\tau'(q)=u(\bm{z}(q))$, as claimed.

    Now, \cref{im:decr} follows from \cref{im:deriv} since $\tau$ is concave and therefore has monotonically decreasing derivative.

    To see \cref{im:fa}, just as in the proof of \cref{p:conc-char}, fix $\alpha\in[u(\bm{z}(\infty)),u(\bm{z}(-\infty))]$ and let $q_\alpha\in\R\cup\{-\infty,\infty\}$ be such that $u(\bm{z}(q_\alpha))=\alpha$.
    Then
    \begin{equation*}
        f(\alpha)\geq v(\bm{z}(q_\alpha))=q_\alpha \alpha-\tau(q)=\tau^*(\alpha)
    \end{equation*}
    since $\alpha\in\partial\tau(q_\alpha)$.
    We recall that the inequality $f(\alpha)\leq\tau^*(\alpha)$ holds in general from \cref{p:conc-char}.

    Finally, \cref{im:outside-dom} follows since
    \begin{equation*}
        u(\bm{z}(\pm\infty))=\lim_{q\to\pm\infty}\frac{\tau(q)}{q}
    \end{equation*}
    so for any $\alpha\notin[u(\bm{z}(\infty)),u(\bm{z}(-\infty))]$, $f(\alpha)\leq\tau^*(\alpha)=-\infty$.
\end{proof}
As a quick example, one setting in which such a function $\bm{z}$ must necessarily exist occurs the optimization defining $\tau(q)$ has a unique minimum for all $q\in\R$.
\begin{corollary}
    Suppose for each $q\in\R$, the minimization defining $\tau(q)$ is attained uniquely at a value $\bm{z}(q)$.
    Then $\bm{z}\colon\R\to\Delta$ is continuous.
\end{corollary}
\begin{proof}
    Since $\Delta$ is compact, we may use the sequential characterization of continuity along the sequences for which the limit $\bm{z}(q_n)$ exists.
    Let $q_0\in\R$ be arbitrary and let $(q_n)_{n=1}^\infty$ be any sequence of real numbers converging to $q_0$ such that the limit $\bm{w}=\lim_{n\to\infty}\bm{z}(q_n)$ exists.
    Then since $\tau$ is a continuous, $u$ is continuous, and $v$ is upper semicontinuous,
    \begin{equation*}
        \tau(q_0)\leq q_0\cdot u(\bm{w})-v(\bm{w})\leq \lim_{n\to\infty}(q_n\cdot u(\bm{z}(q_n))-v(\bm{z}(q_n)))=\lim_{n\to\infty}\tau(q_n)=\tau(q_0)
    \end{equation*}
    so that all the inequalities are in fact equalities.
    But $\bm{z}(q_0)$ is unique with this property, so $\bm{z}(q_0)=\bm{w}$ and continuity follows.
\end{proof}
\begin{acknowledgements}
    The author thanks Kenneth Falconer for many comments on various draft versions of this paper, which have helped to improve the exposition greatly.
    He also thanks Amlan Banaji for pointing out various typos and suggesting some improvements, as well as the anonymous referee for a careful reading which caught many small errors and inconsistencies.
    This paper has benefited greatly from conversations with István Kolossváry (who introduced the author to the method of types) and with Thomas Jordan.
\end{acknowledgements}
\end{document}

%% file: figures/lq_plot_asymp.tex
\begin{tikzpicture}[>=stealth,yscale=1,xscale=2.2]
    \draw[->] (-2,0) -- (4,0);
    \draw[->] (0,-4) -- (0,3);

    \draw[thick] plot file {figures/lq_points.txt};
    \draw[thick, gray, dashed] (0,0) -- (-2,-4);
    \draw[thick, gray, dashed] (0,-0.438018) -- (4,4*0.63093-0.438018);

    \node[right] at (-1.8, -4) {$q\mapsto q\cdot \kappa_{\max} -t_{\max}$};
    \node[left] at (3.6,4*0.63093-0.438018) {$q\mapsto q\cdot \kappa_{\min} -t_{\min}$};

    \node[above left, fill=white] at (0,0) {$-t_{\max}$};
    \draw[dotted] (0,-0.438018) -- (1,-0.438018) node[right] {$-t_{\min}$} ;
    \draw[dotted] (1,0) -- (1,0.5) node[above] {$(1,0)$};
    \node[vtx] at (1,0) {};
    \node[vtx] at (0,-0.943936) {};
    \node[below right] at (0, -0.943936) {$-\dimH K$};

    \node[vtx] at (0,0) {};
    \node[vtx] at (0,-0.438018) {};
\end{tikzpicture}

%% file: figures/multifractal.tex
\begin{tikzpicture}[scale=4,>=stealth]
    \draw[->] (-0.1,0) -- (2.2,0);
    \draw[->] (0,-0.1) -- (0,1.1);
    \draw[thick] plot file {figures/multifractal_points.txt};

    \coordinate (lp) at (0.6309297535714573, 0.4380178794859424);
    \coordinate (rp) at (2,0);
    \coordinate (midpoint) at (0.817455,0.817455);

    \node[vtx] at (lp) {};
    \node[vtx] at (rp) {};
    \node[vtx] at (midpoint) {};

    \node[below] at (lp) {$(\kappa_{\min},t_{\min})$};
    \node[above right] at (rp) {$(\kappa_{\max},t_{\max})$};

    \draw[dashed] (1.6,0.943936) -- (0,0.943936) node[left]{$\dimH K$};
    \draw[dotted,font=\small] (-0.05,-0.05) -- (1.1,1.1) node[above right]{$y=x$};
    \draw[dashed] (0.817455,0.817455) -- (0,0.817455) node[left]{$\dimH\mu_{\bm{p}}$};

    \draw[dashed] (1.11175,0.943936) -- (1.11175,0) node[below]{$\tau_{\bm{p}}'(0)$};
\end{tikzpicture}

%% file: figures/simplex.tex
\begin{tikzpicture}[scale=4]
    \draw[thick, gray] (-0.5,-0.8660254) -- (-0.5, 0.8660254) -- (1,0) -- cycle;
    \draw[thick] plot file {figures/simplex_points.txt};
    \draw[dotted] plot file {figures/contour_points_2.txt};
    \draw[dotted] plot file {figures/contour_points_3.txt};
    \draw[dotted] plot file {figures/contour_points_4.txt};

    \coordinate (zminusinf) at (-0.5, -0.866025);
    \coordinate (zplusinf) at (0.427051, 0.330792);

    \node[vtx] at (zminusinf) {};
    \node[below left] at (zminusinf) {$\bm{z}(-\infty)$};
    \node[above right] at (zplusinf) {$\bm{z}(\infty)$};

    \node[vtx] at (zplusinf) {};

    \node[vtx] at (0.25, 0.){};
    \node[right] at (0.25, 0){$\bm{z}(q)$};
\end{tikzpicture}

%% file: main.bib
@article{zbl:0929.28007,
  author = {Lau, Ka-Sing and Ngai, Sze-Man},
  doi = {10.1006/aima.1998.1773},
  eprint = {0929.28007},
  eprinttype = {zbl},
  journal = {Adv. Math.},
  language = {English},
  pages = {45--96},
  title = {Multifractal measures and a weak separation condition},
  volume = {141},
  year = {1999},
  zbl = {0929.28007},
  zbmath = {01251308},
}

@article{zbl:1040.28014,
  author = {Olsen, Lars and Winter, Steffen},
  doi = {10.1112/S0024610702003630},
  eprint = {1040.28014},
  eprinttype = {zbl},
  journal = {J. Lond. Math. Soc., II. Ser.},
  language = {English},
  pages = {103--122},
  title = {Normal and non-normal points of self-similar sets and divergence points of self-similar measures},
  volume = {67},
  year = {2003},
  zbl = {1040.28014},
  zbmath = {02072301},
}

@article{zbl:0009.05301,
  author = {Besicovitch, A. S.},
  doi = {10.1112/jlms/s1-9.2.126},
  eprint = {0009.05301},
  eprinttype = {zbl},
  journal = {J. Lond. Math. Soc.},
  language = {English},
  pages = {126--131},
  subtitle = {On rational approximation to real numbers},
  title = {Sets of fractional dimensions {IV}},
  volume = {9},
  year = {1934},
  zbl = {0009.05301},
  zbmath = {03012906},
}

@article{zbl:0014.25702,
  author = {Knichal, Vladimír},
  doi = {10.21136/CPMF.1936.109344},
  eprint = {0014.25702},
  eprinttype = {zbl},
  journal = {Čas. Mat. Fys.},
  language = {German},
  origtitle = {Dyadic expansions and Hausdorff measure},
  pages = {195--209},
  title = {Dyadische Entwicklungen und Hausdorffsches Mass},
  volume = {65},
  year = {1936},
  zbl = {0014.25702},
  zbmath = {03022461},
}

@article{zbl:0031.20801,
  author = {Eggleston, H. G.},
  doi = {10.1093/qmath/os-20.1.31},
  eprint = {0031.20801},
  eprinttype = {zbl},
  journal = {Q. J. Math., Oxf. Ser.},
  language = {English},
  pages = {31--36},
  publisher = {Oxford University Press, Oxford},
  title = {The fractional dimension of a set defined by decimal properties},
  volume = {20},
  year = {1949},
  zbl = {0031.20801},
  zbmath = {03048554},
}

@book{zbl:0193.18401,
  author = {Rockafellar, R. Tyrrell},
  doi = {10.1515/9781400873173},
  edition = {1},
  eprint = {0193.18401},
  eprinttype = {zbl},
  journal = {Princeton Math. Ser.},
  language = {English},
  publisher = {Princeton University Press, Princeton, NJ},
  title = {Convex analysis},
  volume = {28},
  year = {1970},
  zbl = {0193.18401},
  zbmath = {03307135},
}

@article{zbl:0763.58018,
  author = {Cawley, Robert and Mauldin, R. Daniel},
  doi = {10.1016/0001-8708(92)90064-R},
  eprint = {0763.58018},
  eprinttype = {zbl},
  journal = {Adv. Math.},
  language = {English},
  pages = {196--236},
  title = {Multifractal decompositions of {Moran} fractals},
  volume = {92},
  year = {1992},
  zbl = {0763.58018},
  zbmath = {00036405},
}

@book{zbl:0475.28009,
  author = {Walters, Peter},
  eprint = {0475.28009},
  eprinttype = {zbl},
  language = {English},
  publisher = {Springer, Cham},
  series = {Grad. Texts Math.},
  title = {An introduction to ergodic theory},
  volume = {79},
  year = {1982},
  zbl = {0475.28009},
  zbmath = {03745547},
}

@book{zbl:0587.28004,
  address = {Cambridge},
  author = {Falconer, Kenneth J.},
  eprint = {0587.28004},
  eprinttype = {zbl},
  language = {English},
  publisher = {Cambridge University Press},
  series = {Camb. Tracts Math.},
  title = {The geometry of fractal sets},
  volume = {85},
  year = {1985},
  zbl = {0587.28004},
  zbmath = {03941840},
}

@article{zbl:0605.58028,
  author = {Ledrappier, F. and Young, L.-S.},
  doi = {10.2307/1971328},
  eprint = {0605.58028},
  eprinttype = {zbl},
  journal = {Ann. Math.},
  language = {English},
  pages = {509--539},
  subtitle = {Part {I}: {Characterization} of measures satisfying {Pesin}'s entropy formula},
  title = {The metric entropy of diffeomorphisms},
  volume = {122},
  year = {1985},
  zbl = {0605.58028},
  zbmath = {03977985},
}

@article{zbl:0664.58022,
  address = {Cambridge},
  author = {Rand, David A.},
  doi = {10.1017/S0143385700005162},
  eprint = {0664.58022},
  eprinttype = {zbl},
  journal = {Ergodic Theory Dyn. Syst.},
  language = {English},
  pages = {527--541},
  publisher = {Cambridge University Press},
  title = {The singularity spectrum {\(f(\alpha)\)} for cookie-cutters},
  volume = {9},
  year = {1989},
  zbl = {0664.58022},
  zbmath = {04086585},
}

@book{zbl:1140.94001,
  address = {Hoboken, NJ},
  author = {Cover, Thomas M. and Thomas, Joy A.},
  doi = {10.1002/047174882X},
  eprint = {1140.94001},
  eprinttype = {zbl},
  language = {English},
  publisher = {John Wiley \& Sons},
  title = {Elements of information theory},
  year = {2006},
  zbl = {1140.94001},
  zbmath = {05060438},
}

@article{zbl:1141.28006,
  author = {Baek, I. S. and Olsen, L. and Snigireva, N.},
  doi = {10.1016/j.aim.2007.02.003},
  eprint = {1141.28006},
  eprinttype = {zbl},
  journal = {Adv. Math.},
  language = {English},
  pages = {267--287},
  title = {Divergence points of self-similar measures and packing dimension},
  volume = {214},
  year = {2007},
  zbl = {1141.28006},
  zbmath = {05180272},
}

@book{zbl:1177.60035,
  author = {Dembo, Amir and Zeitouni, Ofer},
  doi = {10.1007/978-3-642-03311-7},
  eprint = {1177.60035},
  eprinttype = {zbl},
  language = {English},
  publisher = {Berlin: Springer},
  series = {Stoch. Model. Appl. Probab.},
  title = {Large deviations techniques and applications},
  volume = {38},
  year = {2010},
  zbl = {1177.60035},
  zbmath = {05639438},
}

@article{zbl:1184.37028,
  author = {Halsey, Thomas C. and Jensen, Mogens H. and Kadanoff, Leo P. and Procaccia, Itamar and Shraiman, Boris I.},
  doi = {10.1103/PhysRevA.33.1141},
  eprint = {1184.37028},
  eprinttype = {zbl},
  journal = {Phys. Rev. A},
  language = {English},
  pages = {1141--1151},
  title = {Fractal measures and their singularities: the characterization of strange sets},
  volume = {33},
  year = {1986},
  zbl = {1184.37028},
  zbmath = {05695272},
}

@article{zbl:1228.37025,
  author = {Barral, Julien and Feng, De-Jun},
  doi = {10.1088/0951-7715/24/9/010},
  eprint = {1228.37025},
  eprinttype = {zbl},
  journal = {Nonlinearity},
  language = {English},
  pages = {2563--2567},
  title = {Non-uniqueness of ergodic measures with full {Hausdorff} dimension on a {Gatzouras}--{Lalley} carpet},
  volume = {24},
  year = {2011},
  zbl = {1228.37025},
  zbmath = {05956581},
}

@article{zbl:1319.37016,
  arxiv = {1102.5427},
  author = {Climenhaga, Vaughn},
  doi = {10.1017/etds.2014.12},
  eprint = {1319.37016},
  eprinttype = {zbl},
  journal = {Ergodic Theory Dyn. Syst.},
  language = {English},
  pages = {1409--1450},
  publisher = {Cambridge University Press, Cambridge},
  title = {The thermodynamic approach to multifractal analysis},
  volume = {34},
  year = {2014},
  zbl = {1319.37016},
  zbmath = {06366525},
}

@book{zbl:1390.28012,
  address = {Cambridge},
  author = {Bishop, Christopher J. and Peres, Yuval},
  doi = {10.1017/9781316460238},
  eprint = {1390.28012},
  eprinttype = {zbl},
  language = {English},
  number = {162},
  publisher = {Cambridge University Press},
  series = {Camb. Stud. Adv. Math.},
  title = {Fractals in probability and analysis},
  year = {2017},
  zbl = {1390.28012},
  zbmath = {06653783},
}

@article{zbl:1371.37012,
  author = {Ledrappier, F. and Young, L.-S.},
  doi = {10.2307/1971329},
  eprint = {1371.37012},
  eprinttype = {zbl},
  journal = {Ann. Math.},
  language = {English},
  pages = {540--574},
  subtitle = {Part {II}: {Relations} between entropy, exponents and dimension},
  title = {The metric entropy of diffeomorphisms},
  volume = {122},
  year = {1985},
  zbl = {1371.37012},
  zbmath = {06781245},
}

@article{zbl:1387.37026,
  author = {Das, Tushar and Simmons, David},
  doi = {10.1007/s00222-017-0725-5},
  eprint = {1387.37026},
  eprinttype = {zbl},
  journal = {Invent. Math.},
  language = {English},
  pages = {85--134},
  title = {The {Hausdorff} and dynamical dimensions of self-affine sponges: a dimension gap result},
  volume = {210},
  year = {2017},
  zbl = {1387.37026},
  zbmath = {06798848},
}

@article{zbl:1426.11079,
  author = {Shmerkin, Pablo},
  doi = {10.4007/annals.2019.189.2.1},
  eprint = {1426.11079},
  eprinttype = {zbl},
  journal = {Ann. Math.},
  language = {English},
  pages = {319--391},
  title = {On {Furstenberg}'s intersection conjecture, self-similar measures, and the {{\(L^q\)}} norms of convolutions},
  volume = {189},
  year = {2019},
  zbl = {1426.11079},
  zbmath = {07041748},
}

@book{zbl:0819.28004,
  author = {Mattila, Pertti},
  eprint = {0819.28004},
  eprinttype = {zbl},
  language = {English},
  publisher = {Cambridge: Univ. Press},
  series = {Camb. Stud. Adv. Math.},
  subtitle = {Fractals and rectifiability},
  title = {Geometry of sets and measures in {Euclidean} spaces},
  volume = {44},
  year = {1995},
  zbl = {0819.28004},
  zbmath = {00739280},
}

@article{zbl:1549.37013,
  address = {London},
  arxiv = {2205.01043},
  author = {Kolossváry, István},
  doi = {10.1112/jlms.12767},
  eprint = {1549.37013},
  eprinttype = {zbl},
  journal = {J. Lond. Math. Soc.},
  language = {English},
  pages = {666--701},
  publisher = {London Mathematical Society},
  title = {The \(L^q\) spectrum of self-affine measures on sponges},
  volume = {108},
  year = {2023},
  zbl = {1549.37013},
  zbmath = {07731132},
}

@book{zbl:1543.28001,
  address = {Providence, RI},
  author = {Bárány, Balázs and Simon, Károly and Solomyak, Boris},
  doi = {10.1090/surv/276},
  eprint = {1543.28001},
  eprinttype = {zbl},
  language = {English},
  publisher = {American Mathematical Society},
  series = {Math. Surv. Monogr.},
  title = {Self-similar and self-affine sets and measures},
  volume = {276},
  year = {2023},
  zbl = {1543.28001},
  zbmath = {07759243},
}

@article{zbl:0873.28003,
  author = {Arbeiter, Matthias and Patzschke, Norbert},
  doi = {10.1002/mana.3211810102},
  eprint = {0873.28003},
  eprinttype = {zbl},
  journal = {Math. Nachr.},
  language = {English},
  pages = {5--42},
  title = {Random self-similar multifractals},
  volume = {181},
  year = {1996},
  zbl = {0873.28003},
  zbmath = {00927263},
}

@book{zbl:0869.28003,
  address = {Chichester},
  author = {Falconer, Kenneth},
  eprint = {0869.28003},
  eprinttype = {zbl},
  language = {English},
  publisher = {John Wiley \& Sons},
  title = {Techniques in fractal geometry},
  year = {1997},
  zbl = {0869.28003},
  zbmath = {00981057},
}
